\documentclass[11pt]{amsart}

\usepackage{amscd}
\usepackage{amssymb}

\usepackage[frame,curve,matrix,arrow]{xy}
\CompileMatrices   

\swapnumbers
\usepackage{amsthm}

\newtheorem{thm}[equation]{Theorem}
\newtheorem{lem}[equation]{Lemma}
\newtheorem{cor}[equation]{Corollary}
\newtheorem{prop}[equation]{Proposition}

\newtheorem*{thm*}{Theorem}
\newtheorem*{prop*}{Proposition}
\newtheorem*{cor*}{Corollary}
\newtheorem*{lem*}{Lemma}

\theoremstyle{definition} %

\newtheorem*{defn*}{Definition}

\newtheorem{eg}[equation]{Example}

\theoremstyle{remark} %
\newtheorem{rmk}[equation]{Remark}

\newtheorem*{rmk*}{Remark}
\newtheorem*{rmks*}{Remarks}

\newtheoremstyle{exercise}
  {3pt}
  {3pt}
  {\small}
  {\parindent}
  {\sc\small}
  {.}
  {.5em}
   {}     
  {}

\theoremstyle{exercise}

\numberwithin{equation}{section}

\makeatletter
{\renewcommand{\theequation}{#1}}%
{\renewcommand{\theequation}{\arabic{equation}}\addtocounter{equation}{-1}\global\@ignoretrue}
\makeatother

\makeatletter
{\renewcommand{\theequation}{#1}\begin{eqnarray}}%
{\end{eqnarray}\renewcommand{\theequation}{\arabic{equation}}\addtocounter{equation}{-1}\global\@ignoretrue}
\makeatother

\makeatletter
\newenvironment{borel}[1]%
{\smallskip \refstepcounter{equation}\noindent{\textbf \theequation. }{{\textbf{#1.}}}}%
{\smallskip \global\@ignoretrue}
\makeatother

\makeatletter
{\smallskip \refstepcounter{equation}{\sc \theequation}{\sc (#1).}}%
{\smallskip \global\@ignoretrue}
\makeatother

\makeatletter
{\smallskip \refstepcounter{equation}\noindent{\sc \theequation.}{\sl{ #1.}}}%
{\smallskip \global\@ignoretrue}
\makeatother

\makeatletter
\newenvironment{borel*}%
{\smallskip \refstepcounter{equation}\noindent{\textbf \theequation.}}%
{\global\@ignoretrue}
\makeatother

\newcommand{\flist}[1]{\hangindent\leftmargini\textup{(1)}\hskip\labelsep {#1}%
\begin{enumerate}%
\setcounter{enumi}{1}%
}

\newcommand{\ot}{\otimes}

\newcommand{\C}{{\mathbb{C}}}        

\newcommand{\e}{\varepsilon}

\newcommand{\la}{\lambda}
\newcommand{\D}{\Delta}

\newcommand{\G}{{\Gamma}}       

\newcommand{\oddots}{{\mathinner{\mkern1mu\raise1pt\vbox{\kern7pt\hbox{.}}\mkern2mu\raise4pt\hbox{.}\mkern2mu\raise7pt\hbox{.}\mkern1mu}}}

\newcommand{\s}{\sigma}

\newcommand{\qform}[1]{{\langle{#1}\rangle}}                   

\DeclareMathOperator{\Spin}{Spin}           

\newcommand{\Ga}{\mathbb{G}_a}
\newcommand{\Gm}{\mathbb{G}_m}

\DeclareMathOperator{\im}{im}

\newcommand{\iso}{\xrightarrow{\sim}}

\newcommand{\ra}{\rightarrow}

\newcommand{\wbar}{\overline{\omega}}

\newcommand{\anmn}{{\{ \alpha_{n-1}, \alpha_n \}}}

\newcommand{\hw}{v^+}
\newcommand{\lw}{v^-}

\newcommand{\cP}{\mathcal{P}}

\newcommand{\ach}{\check{\alpha}}

\newcommand{\g}{\mathfrak{g}}

\theoremstyle{remark}
\newtheorem*{brmk}{Bibliographic remarks}

\DeclareMathOperator{\sub}{SubSp}
\DeclareMathOperator{\inid}{Inid}

\newcommand{\atil}{\widetilde{\alpha}}

\newcommand{\Zbar}{\overline{Z}}

\newcommand{\X}{\mathfrak{X}}
\newcommand{\B}{\mathcal{B}}

\newcommand{\kbar}{\overline{k}}

\newcommand{\darkradG}{0.115}

\newcommand{\kspan}{\text{$k$-span}\,}

\begin{document}

\author{Michael Carr}
\address{Department of Mathematics, University of Michigan, Ann Arbor, MI 48109-1109}
\email{mpcarr@umich.edu}

\title%
{Geometries, the principle of duality, and algebraic groups}
\author{Skip Garibaldi}
\thanks{Corresponding author: Garibaldi}
\address{Department of Mathematics \& Computer Science, Emory University, Atlanta, GA 30322}
\email{skip@member.ams.org}
\urladdr{http://www.mathcs.emory.edu/{\textasciitilde}skip/}

\begin{abstract}
J.~Tits gave a general recipe for producing an abstract geometry from a semisimple algebraic group.  This expository paper describes a uniform method for giving a concrete realization of Tits's geometry and works through several examples. We also give a criterion for recognizing the automorphism of the geometry induced by an automorphism of the group.  The $E_6$ geometry is studied in depth.
\end{abstract}

\date{\today}

\maketitle

\setcounter{tocdepth}{1}
\tableofcontents

J.~Tits's theory of buildings associated with semisimple algebraic groups gives a unified method of extracting a geometry from a group.  For example, the group $SL_n$ gives rise to $(n-1)$-dimensional projective space.  Tits's geometry however is very abstract.  Speaking precisely, one obtains an \emph{incidence geometry}, which consists of an abstract set of objects each with a given type, and a reflexive, symmetric binary relation on the set of objects called incidence.  We find it more palatable to think of projective space in a concrete way, as the collection of subspaces of some explicit vector space.  In Section \ref{realize} we give an explicit recipe for concretizing Tits's incidence geometry.  

The midsection of this paper consists of explicit descriptions of the concrete realizations of the geometries for split groups of type $A$, $D$, $E_6$, $E_7$, $F_4$, and $G_2$.  (Readers should have little trouble filling in the missing types $B$ and $C$.  We do not know a good description for the geometry of type $E_8$, but we make a few comments in \ref{E8}.)  Such descriptions may be found in a variety of places in the literature, e.g., \cite{Cohen:pl} or \cite{FF}.  The main innovation here is that our recipe produces a realization of the geometry by a largely deterministic process beginning from the root system of the group and a fundamental representation, whereas approaches in the literature have the appearance of being ad hoc.  Roughly speaking, Tits developed his theory of buildings by abstracting and unifying known properties of the various concrete geometries \cite[p.~v]{Ti:BN}, so our approach here reverses the historical development.

Our principal tool is the representation theory of semisimple groups; we only use the most elementary results, but we exploit those ruthlessly.  Consequently, throughout this paper, \emph{our base field $k$ is assumed to have characteristic zero}.  Some results hold over an arbitrary field, see Remarks \ref{min.free} and \ref{jord.free}.

The final portion of this paper concerns duality.  \'Elie Cartan was already aware (see \cite[p.~362]{Cartan}) that the ``outer'' automorphism $g \mapsto 
(g^{-1})^t$ (where $t$ denotes transposition) of $SL_3$ gives rise to a polarity in the projective plane and so to the principle of duality.  In general, an outer automorphism of a semisimple group gives an automorphism $\psi$ of the corresponding geometry, hence also a ``principle of duality'' (or ``triality'' or ...).  In Sections \ref{outer}--\ref{E6.dual} below, we give a criterion for recognizing such an automorphism $\psi$ and apply our criterion to get an explicit description of $\psi$ in essentially all cases.  The explicit description of $\psi$ for the $E_6$ geometry is new.  Having an explicit formula for $\psi$ is useful for giving a concrete description of the projective homogeneous varieties for groups of ``outer type", see for example \cite[p.~172]{MPW2}.

\smallskip

We hope that readers who are not familiar with, say, exceptional groups will find the presentation here unusually accessible because of the uniform treatment in the common language of representation theory.  Experts will note that the concept of inner ideal occurs naturally in the $E_7$ geometry (\S\ref{E7}) and that we do not need to mention octonions at all in our discussion of triality in \S\ref{D4}.

\smallskip

Our hypothetical reader is moderately familiar with the theory of linear algebraic groups as in \cite{Borel}, 
\cite{Sp:LAG}, or \cite{St} and the classification of irreducible modules via highest weight vectors from \cite[Ch.~VI]{Hum:LA}, \cite[\S14]{FH}, \cite[\S5.1]{GW}, or \cite{Var}.

\section{Tits's geometry $\G_P$} \label{Tits}

In this section, we describe Tits's recipe for producing an incidence geometry from a certain kind of algebraic group.  An \emph{incidence geometry} is a set of objects, each of some type (e.g., point, line), together with a symmetric binary relation known as \emph{incidence}.  There is just one further axiom: objects of the same type are incident if and only if they are equal.

\begin{rmk*} Modern formulations of Tits's recipe take a group and construct a building rather than an incidence geometry.  From our perspective, a building is an incidence geometry with extra structure that we do not really need.  So we deal only with the much-simpler-to-define incidence geometries as in \cite{Ti:ex}.  For a presentation in terms of buildings, see \cite{Ti:BN}, \cite[42.3.6]{TW}, \cite[Ch.~V]{Brown}, or \cite[\S4.4]{Scharlau}.
\end{rmk*}

We start with a root system $\Phi$ in the sense of \cite[\S1]{St}, a ``reduced root system'' in the language of \cite{Bou:g4}.  There is a ``split" simply connected algebraic group $G$ with root system $\Phi$, and it is uniquely determined up to isomorphism. Taking $\Phi$ of type $A_n$, we obtain a group $G$ isomorphic to $SL_{n+1}$.  
When $k$ is algebraically closed, every semisimple 
algebraic group is obtained in this fashion, or is a quotient of a group obtained in this fashion.  For example, the group $SO_{2n}(\C)$ is a quotient of $\Spin_{2n}(\C)$, which is constructed from the root system $D_n$.  

\borel{Parabolic subgroups}  A \emph{Borel subgroup} is a maximal closed, connected, solvable subgroup of $G$.  Fix one and call it $B$; it (combined with a split maximal torus $T$ contained in it) determines a set of simple roots $\D$ in $\Phi$.  We abuse notation by writing $\D$ also for the associated Dynkin diagram.

A closed subgroup of $G$ is called \emph{parabolic} if it contains a Borel subgroup, and we call a parabolic subgroup \emph{standard} if it contains the Borel $B$.  There is an inclusion-reversing bijection
   \[
   \begin{array}{r@{\quad \leftrightarrow\quad}l}
   \text{\fbox{standard parabolic subgroups of $G$}}&\text{\fbox{subsets of $\D$}} \\
   B&\D \\
   G&\emptyset
   \end{array}
   \]
The maximal proper standard parabolics are in one-to-one correspondence with the elements of $\D$.  We write $P_\delta$ for the standard parabolic corresponding to $\delta \in \D$.

\begin{eg}[Parabolics in $SL_4$] \label{SL4.1}
As an illustration, the Dynkin diagram of $SL_4$ is 
   \[
   \begin{picture}(4,0.7)
    \put(1,0.5){\line(3,0){2}}

    \multiput(1,0.5)(1,0){3}{\circle*{\darkradG}}

    \put(1,0.1){\makebox(0,0.4)[b]{$\alpha_1$}}
    \put(2,0.1){\makebox(0,0.4)[b]{$\alpha_2$}}
    \put(3,0.1){\makebox(0,0.4)[b]{$\alpha_3$}}
\end{picture}
   \]
Here we have labeled the vertices as in \cite{Bou:g4}.  The upper triangular matrices are a Borel subgroup, which we take to be $B$.  For our maximal torus $T$, we take the diagonal matrices.  The maximal proper standard parabolics have the form
   \begin{equation} \label{SL4.para}
   P_{\alpha_1} = \left( \begin{smallmatrix}
   *&*&*&* \\
   0&*&*&* \\
   0&*&*&* \\
   0&*&*&* 
   \end{smallmatrix} \right), \quad
  P_{\alpha_2} = \left( \begin{smallmatrix}
   *&*&*&* \\
   *&*&*&* \\
   0&0&*&* \\
   0&0&*&* 
   \end{smallmatrix} \right), \quad
  P_{\alpha_3} = \left( \begin{smallmatrix}
   *&*&*&* \\
   *&*&*&* \\
   *&*&*&* \\
   0&0&0&* 
   \end{smallmatrix} \right).
   \end{equation}
\end{eg}

\borel{Tits's geometry $\G_P$  \textrm{(\cite[\S5]{Ti:BN}, \cite[42.3.6]{TW})}} \label{Tits.geom}
Tits defines the objects of the incidence geometry to be the maximal proper parabolic subgroups of $G$. 
Since the Borel subgroups of $G$ are all conjugate, every parabolic is conjugate to a unique standard parabolic \cite[21.12]{Borel}.  Therefore, every maximal proper parabolic subgroup corresponds to a unique element $\delta \in \D$; this is the \emph{type} of the parabolic.  Two maximal proper parabolics are said to be \emph{incident} if their intersection contains a parabolic subgroup.  We write $\G_P$ for this incidence geometry, where the subscript $P$ is meant to remind the reader that the objects are parabolic subgroups of $G$.

In the case of $SL_4$, one typically calls the parabolics of type $\alpha_1$ ``points'', $\alpha_2$ ``lines'', and $\alpha_3$ ``planes''.  This identifies the geometry $\G_P$ with 3-dimensional projective space.

Note that the notion of incidence we defined for $\G_P$ satisfies the axiom from the beginning of the section: 
If parabolics $P$ and $P'$ of the same type are incident, then $P \cap P'$ contains a Borel subgroup, which after conjugation we may assume is the standard Borel $B$.  But by hypothesis, $P$ and $P'$ have the same type, hence $P = P'$.

\section{A concrete geometry $\G_V$, part I} \label{realize}

We now take Tits's incidence geometry $\G_P$---whose objects are certain parabolic subgroups of $G$---and produce an isomorphic incidence geometry $\G_V$ whose objects are subspaces of a fixed vector space $V$.  For example, in the case where $G$ is $SL_n$, $\G_V$ consists of the nonzero, proper subspaces of $k^n$.

Fix a representation of $G$ on a finite-dimensional vector space $V$, i.e., a homomorphism of algebraic groups $G \ra GL(V)$.  We assume that the representation is nontrivial (i.e., the image of $G$ is not just the identity transformation) and irreducible (i.e., there is no nonzero, proper $G$-invariant subspace).  For each vertex $\delta$ of the Dynkin diagram, we choose a nonzero, proper subspace $V_\delta$ of $V$ that is invariant under the parabolic $P_\delta$.\footnote{For the moment, we assume that such a subspace exists.  The doubting reader may wish to glance ahead at Prop.~\ref{stable}.}
Every maximal proper parabolic subgroup $P'$ is conjugate to a unique standard parabolic subgroup $P_\delta$, and we define a subspace $V_{P'}$ via
   \[
   V_{P'} := g V_\delta \quad \text{for $g \in G(k)$ such that $P' = g P_\delta g^{-1}$.}
   \]
Of course, $g$ is not uniquely determined, but if $h\in G(k)$ also satisfies $h P_\delta h^{-1} = P'$, then $h^{-1} g$ normalizes $P_\delta$.  Since $P_\delta$ is its own normalizer \cite[11.16]{Borel},  $g$ equals $hp$ for some $p \in P_\delta(k)$ and $h V_\delta = g V_\delta$.   We remark that the stabilizer of $V_{P'}$ in $G$ is precisely $P'$.  Indeed, the stabilizer of $V_{P'}$ is a closed subgroup of $G$ containing $P'$, hence must be $P'$ or $G$.  Since the representation $V$ is irreducible, the stabilizer is $P'$.

We define an incidence geometry $\G_V$ whose objects are the subspaces $V_P$ of $V$, as $P$ ranges over the maximal proper parabolic subgroups of $G$.  The map $P \mapsto V_P$ is surjective by definition, but it is also injective because $P$ is precisely the stabilizer of $V_P$.  Therefore, the sets of objects in $\G_P$ and $\G_V$ are isomorphic.  Define the notions of type and incidence in $\G_V$ by transporting them from $\G_P$.  Speaking precisely, we say that $V_P$ in $\G_V$ has type $\delta \in \D$ if the parabolic $P$ is of type $\delta$.  We define two objects in $\G_V$ to be incident if and only if the corresponding parabolic subgroups are incident. 

In this very simple way, we have obtained a realization of Tits's abstract geometry $\G_P$ as a collection of subspaces of the concrete vector space $V$.  This recipe begs two obvious questions:
\begin{align}
&\parbox{4in}{Are we guaranteed that a nonzero, proper $P_\delta$-invariant subspace of $V$ exists?}  \label{ques.1}\\
&\parbox{4in}{Is there a way to tell if two subspaces of $V$ in $\G_V$ are incident without discussing the corresponding parabolic subgroups?}\label{ques.2}
\end{align}
We will address these two questions in the next section and the examples that follow it.  But first, here is an example to illustrate the construction.

\begin{eg}[$\G_V$ for $SL_4$] \label{SL4}
Referring to the description of the parabolic subgroups of $SL_4$ from Example \ref{SL4.1}, we see that $P_{\alpha_i}$ stabilizes the $i$-dimensional subspace of $k^4$ consisting of vectors whose only nonzero entries are in the first $i$ coordinates.  We can take this to be $V_{\alpha_i}$.  The objects of $\G_V$ are all proper, nonzero subspaces of $V$.

We claim that two elements of $\G_V$ are incident if and only if one is contained in the other. Indeed, let $W$, $W'$ be proper, nonzero subspaces, stabilized by maximal proper parabolics $P$, $P'$ in $SL_4$.  If $W$ and $W'$ are incident, then $P \cap P'$ contains a Borel subgroup.  After conjugation we may assume that $P$ and $P'$ are standard, hence appear in \eqref{SL4.para}.  Then clearly $W$ is contained in $W'$ or vice versa.  Conversely, if $W$ is contained in $W'$, there is some  $g \in SL_4(k)$ such that
$gW$, $gW'$ are equal to $V_{\alpha_i}$, $V_{\alpha_{i'}}$ for some $i \le i'$.  Then $gPg^{-1}$ and $gP'g^{-1}$ equal $P_{\alpha_i}$ and $P_{\alpha_{i'}}$ respectively, hence the parabolics $P$, $P'$ are incident.
\end{eg}

\section{A concrete geometry $\G_V$, part II} \label{realize.2}

We will now make the geometry $\G_V$ from the previous section more concrete by focusing on the case where $V$ is a fundamental irreducible representation of $G$ (i.e., the highest weight of $V$ is a fundamental weight).  We completely answer \eqref{ques.1} in the affirmative in Prop.~\ref{stable} and we partially answer \eqref{ques.2} in Prop.~\ref{incident.1}.

We view $G$ as being constructed from the root system $\Phi$ by the Chevalley construction as in \cite[\S6]{St}.  That is, it is generated by the images of homomorphisms $x_\alpha \!: \Ga \ra G$  as $\alpha$ ranges over the roots in $\Phi$.  Write $\X_\alpha$ for the image of $x_\alpha$.  For each $\alpha \in \Phi$, the map $t \mapsto x_\alpha(t) x_{-\alpha}(-t^{-1}) x_\alpha(1 - t) x_{-\alpha}(1) x_\alpha(-1)$ defines a homomorphism $\Gm \ra G$, which we denote by $h_\alpha$.  The images of the $h_\alpha$'s generate a maximal torus $T$ in $G$.  We fix a set of simple roots $\D$ in $\Phi$ and choose our standard Borel subgroup $B$ to be the one generated by $T$ and the $\X_\rho$ for $\rho$ a positive root.

Fix a root $\beta \in \D$ and let $\omega$ be the corresponding fundamental weight.  In this section, $V$ denotes the irreducible representation with highest weight $\omega$ with respect to our choice of torus $T$ and Borel $B$.  Fix a highest weight vector $\hw$ in $V$.

\begin{borel*} 
For a simple root $\delta \in \D \setminus \{ \beta \}$, we define the \emph{$\delta$-component} of $\D$ to be the connected component of $\beta$ in $\D \setminus \{ \delta \}$.  The \emph{$\beta$-component} is the empty set.

Now fix a $\delta \in \D$.  For each root $\rho$ of $G$, write $\rho = \sum_{\alpha \in \D} c_\alpha \alpha$ for integers $c_\alpha$; we define the \emph{$\delta$-height} of $\rho$ to be $\sum c_\alpha$ as $\alpha$ ranges over the simple roots \emph{not} in the $\delta$-component.  In the case $\delta = \beta$, this is the usual notion of height.

Write $L_\delta$ for the subgroup of $G$ generated by the root subgroups $\X_\rho$ as $\rho$ ranges over the roots of $\delta$-height zero.  The description of $G$ in terms of generators and relations shows that $L_\delta$ is a simple group whose Dynkin diagram is the $\delta$-component.  (It is the semisimple part of the Levi subgroup of the parabolic corresponding to the complement of the $\delta$-component.)  We set $V_\delta$ to be the subspace $L_\delta \, \hw$ spanned by the $L_\delta$-orbit of $\hw$.  We remark that $L_\beta$ is the group with one element and $V_\beta$ is the line $k \hw$.
\end{borel*}

\begin{lem} \label{not.comp}
For each $\delta \in \D$, the subspace $V_\delta$ is a direct sum of the weight spaces in $V$ with weights of the form $\omega - \alpha$ where $\alpha$ has $\delta$-height zero.
\end{lem}

Recall that every weight of $V$ is of the form $\omega - \alpha$ for $\alpha$ a sum of positive roots.
For the proof, we will use repeatedly \cite[p.~209, Lemma 72]{St}, which says: For $v \in V$ of weight $\omega - \alpha$ and $\rho$ a root,
\begin{equation} \label{stlem}
x_{-\rho}(t) v = v + \sum_{i \ge 1} t^i v_i,
\end{equation}
for some $v_i \in V$ with weight $\omega - (\alpha + i \rho)$.  

\begin{proof} 
Write $V$ as $V_0 \oplus V_1$ where $V_0$ (resp., $V_1$) is spanned by those weight vectors with weight $\omega - \gamma$ where $\gamma$ has $\delta$-height zero (resp., positive $\delta$-height).  We want to show that $V_0$ equals $V_\delta$.

Similarly, write $U_0$ (resp., $U_1$) for the subgroup of $G$ generated by $\X_{-\rho}$, as $\rho$ varies over the positive roots with $\delta$-height zero (resp., positive $\delta$-height).  Equation \eqref{stlem} shows that $V_0$ is $U_0$-invariant, $V_1$ is $U_1$-invariant, and the induced action of $U_1$ on $V/V_1$ is trivial.

Let $U$ denote the subgroup of $G$ generated by the $\X_{-\rho}$ as $\rho$ varies over all the positive roots.  Take a $u \in U(k)$ and write $u = u_1 u_0$ where $u_i$ belongs to $U_i(k)$; this is possible by \cite[21.9]{Borel}.  We have
   \[
   u \hw = u_1 u_0 \hw = u_0 \hw + v_1 \quad \text{for some $v_1 \in V_1$.}
   \]

Because $UB$ is dense in $G$ \cite[14.14]{Borel}, every linear function on the $G$-orbit of $\hw$ also vanishes on the $UB$-orbit of $\hw$.  It follows that the subspace $UB \hw$ is all of $V$.  Since the line $k\hw$ is $B$-invariant, we have observed the standard fact that $U \hw$ is all of $V$.  In particular, the $u \hw$ from the preceding paragraph span $V$ as $u$ ranges over $U(k)$, and we have proved that $V_0$ equals $U_0 \hw$.

The same argument---with $G$, $U$, $B$ replaced with $L_\delta$, $U_0$, and the subgroup of $G$ generated by $T \cap L_\delta$ and the $\X_\rho$ for $\rho$ a positive root of $\delta$-height zero---shows that $V_\delta$ equals $U_0 \hw$.  Hence $V_\delta$ equals $V_0$.
\end{proof}

We now address Question \eqref{ques.1}.

\begin{prop} \label{stable}
For every $\delta \in \D$, the subspace $V_\delta$ is a nonzero, proper subspace of $V$ stabilized by $P_\delta$.
\end{prop}

\begin{proof}
$V_\delta$ is clearly nonzero because it contains $\hw$.  We now show that it is proper.
The highest weight $\atil$ of $\Phi$ equals $\sum_{\alpha \in \D} c_\alpha \alpha$ with every $c_\alpha$ a natural number \cite[10.4A]{Hum:LA}.  In particular, $\qform{\omega, \atil} = c_\beta$ is a positive integer.  Since the set of weights of $V$ is saturated  \cite[\S{VIII.7.2}]{Bou:g4}, $V$ contains a nonzero vector of weight $\omega - \atil$.  But $\atil$  has positive $\delta$-height, so $\omega - \atil$ is not a weight of $V_\delta$.

Finally we show that $V_\delta$ is stabilized by $P_\delta$.  By \cite[14.18]{Borel}, $P_\delta$ is generated by four types of subgroups:
\begin{enumerate}
\item the torus $T$.  It normalizes $L_\delta$ and the line $k\hw$, hence also $V_\delta$.

\item groups $\X_\rho$, for $\rho$ a root with $\delta$-height zero.  These groups belong to $L_\delta$ and so stabilize $V_\delta$.

\item groups $\X_\rho$, for $\rho$ a positive root with nonzero (hence positive) $\delta$-height.  For $\alpha$ a sum of positive roots (possibly zero) with $\delta$-height zero, the sum $\alpha - i \rho$ has negative $\delta$-height for all $i \ge 1$, hence there are no nonzero vectors in $V$ with weight $\omega - (\alpha - i \rho)$.  Equation \eqref{stlem} gives that $\X_\rho$ fixes $V_\delta$ elementwise.

\item groups $\X_{-\rho}$ for $\rho$ a root with positive $\delta$-height and of the form $\sum_{\gamma \in \D \setminus \{ \delta \}} c_\gamma \gamma$.  Since $\rho$ is a root, the subset of $\Delta$ on which the coefficients $c$ are nonzero is connected \cite[\S{VI.1.6}]{Bou:g4}, and since $\rho$'s $\delta$-height is positive, we find that $\rho$ is a sum of simple roots $\gamma$ that are not connected to the $\delta$-component.  Hence $\X_{-\rho}$ commutes with $L_\delta$.  Finally, $c_\beta$ is zero, so $\qform{\omega, -\rho} = 0$ and $\X_{-\rho}$ kills $\hw$, that is, $\X_{-\rho}$ kills $V_\delta$.
\end{enumerate}
\end{proof}

Before treating Question \eqref{ques.2} regarding incidence, we observe that we know a lot about the spaces $V_\delta$.  The case $\delta = \beta$ is trivial, so for the rest of this section we fix a $\delta \in \D$ that is not $\beta$. 
  
\begin{prop} \label{fund}
For $\delta \in \D \setminus \{ \beta \}$, the space $V_\delta$ is a fundamental irreducible representation of $L_\delta$.
\end{prop} 

\begin{proof}
Suppose that $v \in V_\delta$ is fixed by $\X_\rho$ for every positive root $\rho$ with zero $\delta$-height.  The argument in item (3) in the proof of Prop.~\ref{stable} shows that $v$ is fixed by $\X_\rho$ for 
every positive root $\rho$.
Since $V$ is an irreducible representation of $G$, $v$ is in the $k$-span of $\hw$.

The previous paragraph shows that $k \hw$ is the only highest weight line for $V_\delta$ relative to the maximal torus $T_\delta = T \cap L_\delta$ of $L_\delta$.  
Since $V_\delta$ is a completely reducible representation of $L_\delta$, it is irreducible.

The highest weight of $V_\delta$ is the restriction of $\omega$ to $T_\delta$; we denote it by $\wbar$.  Since $\omega$ is a fundamental weight of $G$, the restriction $\wbar$ is a fundamental weight of $L_\delta$.
\end{proof}

The dimension of $V_\delta$ can be looked up in, e.g., \cite[chap.~8, Table 2]{Bou:g4}.

\begin{borel*} \label{spindle}
Recall that the weights are partially ordered by setting $\la_1 \ge \la_2$ if $\la_1 - \la_2$ is a sum of positive weights.  For example, we have already used that $\omega$ is the largest weight of $V$.  On the other hand, $w_0 \omega$ is the smallest weight, where $w_0$ is the longest element of the Weyl group, i.e., the one that sends $\D$ to $-\D$.  (To see this, note that $w_0$ permutes the weights of $V$ and reverses the partial ordering.)

We may apply the same observations to the irreducible representation $V_\delta$ of $L_\delta$. Computing relative to the torus $T_\delta$, the weights of $V_\delta$ lie between the highest weight $\wbar$ and the lowest weight $w_0 \wbar$, where $w_0$ is the longest element in the Weyl group of $L_\delta$.  From the tables in \cite{Bou:g4}, one can quickly find the nonnegative integers $k_\alpha$ such that
   \[
   w_0 \wbar = \wbar - \sum k_\alpha \alpha
   \]
where $\alpha$ runs over the roots in the $\delta$-component.  Considering $V_\delta$ as a subspace of the representation $V$ of $G$, Lemma \ref{not.comp} gives that the weights of $V_\delta$ are precisely those weights $\mu$ of $V$ such that 
    \[
    \omega - \sum k_\alpha \alpha \le \mu \le \omega.
    \]
\end{borel*} 

We close this section by answering the question of incidence---i.e., Question \eqref{ques.2}---in the most important case.  We say that a vertex in a graph is \emph{terminal} if it is joined to at most one other vertex.

We call a subspace $X \in \G_V$ of type $\delta$---i.e., a subspace $X$ in the $G(k)$-orbit of $V_\delta$---a \emph{$\delta$-space}.   (We allow here also the possibility that $\delta = \beta$.)

\begin{prop} \label{incident.1}
Let $X$ be a $\delta$-space in $\G_V$ (with $\delta \ne \beta$).  Suppose that the $\delta$-component of $\D$ is of type $A$ and $\beta$ is a terminal vertex of the $\delta$-component.
\begin{enumerate}
\item Every $1$-dimensional subspace of $X$ is a $\beta$-space.
\item If the $\delta$-component contains the $\delta'$-component, we have: $X$ is incident to a $\delta'$-space $X'$ if and only if $X$ contains $X'$.
\end{enumerate}
\end{prop}

\begin{proof}
To prove the proposition, we may conjugate $X$ and so assume that $X$ is actually $V_\delta$.  The group $L_\delta$ is a special linear group because the $\delta$-component is of type $A$.  Examining the highest weight of the representation $V_\delta$ of $L_\delta$ as in the proof of Prop.~\ref{fund}, we find that there is an isomorphism $L_\delta \iso SL(V_\delta)$ that identifies the natural representations of $L_\delta$ and $SL(V_\delta)$ on $V_\delta$.  Since every 1-dimensional subspace of $V_\delta$ is in the $SL(V_\delta)$-orbit of $V_\beta$, this proves (1).

Now we prove (2).  First suppose that $X$ and $X'$ are incident.  That is, there is some $g \in G(k)$ such that conjugating the stabilizers of $X$ and $X'$ by $g$ gives $P_\delta$ and $P_{\delta'}$ respectively.  By the bijection between parabolics and objects in $\G_V$, we find $gX = V_\delta$ and $gX' = V_{\delta'}$. Since $L_{\delta'}$ is contained in $L_\delta$, clearly $V_{\delta'}$ is contained in $V_\delta$, hence $X'$ is contained in $X$.

Conversely, suppose that $X'$ is contained in $X$.  Applying an element of $G(k)$, we may assume that $X$ equals $V_\delta$. Since $L_\delta$ is $SL(V_\delta)$, all subspaces of $V_\delta$ with the same dimension are $L_\delta$-equivalent.  Hence $X'$ is $L_\delta$-equivalent to $V_{\delta'}$.  The parabolics $P_\delta$, $P_{\delta'}$ corresponding to $V_\delta$, $V_{\delta'}$ contain the standard Borel $B$, hence are incident.
\end{proof}

\section{Example: type $A$ (projective geometry)} \label{A.eg}

In this section, we describe the objects in the geometry $\G_V$ constructed from $G := SL_n$ as in \S\ref{realize.2} when the representation is the standard one on $V := k^n$, corresponding to the simple root $\beta := \alpha_1$.  We ``discover'' that $\G_V$ is projective $(n-1)$-space.  We could do this explicitly in terms of matrices as in Example \ref{SL4}, but such an argument would be hard to generalize to other groups.  Instead, we give an algebraic-group- and representation-theoretic argument.

We defined $V_\beta$, a.k.a.~$V_{\alpha_1}$, to be the 1-dimensional subspace of $V$ spanned by the highest weight vector.  For $\alpha_i \in \D$ with $i \ne 1$,  $V_{\alpha_i}$ is the standard representation of $L_{\alpha_i}$.  Therefore, the dimension of $V_{\alpha_i}$ is precisely $i$.
 We summarize this in the Dynkin diagram, where each vertex is labeled with $\alpha_i$ and $\dim V_{\alpha_i}$:
   \[
   \begin{picture}(9,1)
    \put(1,.5){\line(1,0){3}}
    \put(5.5,.5){\line(1,0){1.5}}
    
    \put(1,0.1){\makebox(0,0.4)[b]{$\beta = \alpha_1$}}  
    \put(2.5,0.1){\makebox(0,0.4)[b]{$\alpha_2$}}
    \put(4,0.1){\makebox(0,0.4)[b]{$\alpha_3$}}
    \put(5.5,0.1){\makebox(0,0.4)[b]{$\alpha_{n-2}$}}
    \put(7,0.1){\makebox(0,0.4)[b]{$\alpha_{n-1}$}}
    \put(1,.6){\makebox(0,0.4)[b]{$1$}}
    \put(2.5,.6){\makebox(0,0.4)[b]{$2$}}
    \put(4,.6){\makebox(0,0.4)[b]{$3$}}
    \put(5.5,.6){\makebox(0,0.4)[b]{${n-2}$}}
    \put(7,.6){\makebox(0,0.4)[b]{${n-1}$}}

    \multiput(4.3,.5)(0.4,0){3}{\circle*{0.015}}

    \multiput(1,.5)(1.5,0){5}{\circle*{\darkradG}}
    
   \end{picture}
   \]
   
Since $SL_n$ acts transitively on the $i$-dimensional subspaces of $V$ for all $i$, we have: the $\alpha_i$-spaces are the $i$-dimensional subspaces of $V$.  By Prop.~\ref{incident.1}.2, two subspaces are incident if and only if one contains the other.   This is the classical description of $(n-1)$-dimensional projective space as consisting of lines through the origin in $k^n$.

\section{Strategy} \label{strategy}

In the next few sections, we will fix a split simply connected group $G$ and give an explicit description of the geometry $\G_V$.
One imagines that the geometry $\G_V$ we have just constructed will be easiest to visualize if the ambient vector space $V$ is small.  With that in mind, we will focus on the case where $V$ is the smallest irreducible representation of $G$.  For $G$ of type $A$, $D_4$, or $E_6$, there are multiple equivalent choices, and we arbitrarily pick one.
\begin{equation} \label{vector}
\begin{array}{c|ccccccccc} \\
\text{type of $G$}&A_n&B_n&C_n&D_n&E_6&E_7&F_4&G_2 \\ \hline
\beta&\alpha_1&\alpha_1&\alpha_1&\alpha_1&\alpha_1&\alpha_7&\alpha_4&\alpha_1 \\
\dim V&n+1&2n+1&2n&2n&27&56&26&7
\end{array}
\end{equation}
We number the elements of $\D$ as in the tables in \cite{Bou:g4}.  We remark that in all cases the root $\beta$ is a terminal vertex of $\D$ as in the hypotheses for Prop.~\ref{incident.1}.  We call a representation $V$ as in the table above a \emph{standard representation} of $G$.  We have omitted type $E_8$, see \ref{E8} for comments.  (We remind the reader that despite our focus on the standard representation of $G$, 
the recipe in \S\ref{realize.2} gives a concrete realization of $\G_P$ for every fundamental representation, and one can compute the dimensions of the $\delta$-spaces using Prop.~\ref{fund}.)

Roughly speaking, each example consists of three parts: dimensions and properties, transitivity, and incidence.

\medskip

In ``dimensions and properties'', for each $\delta \in \D$ we compute the dimension $d$ of $V_\delta$ and some algebraic properties $\cP$ satisfied by $V_\delta$.  These properties will be obviously $G$-invariant, hence they will be satisfied by all the $\delta$-spaces.  Here we restrict ourselves to the tools of elementary representation theory.  This has two advantages.  First, no special background is required to understand the exceptional groups versus the more-familiar classical groups.  Second, we hope the reader will view our descriptions of the $\delta$-spaces as reasonably canonical and not ad hoc.

\medskip

In ``transitivity'', we prove:
\begin{equation} \label{trans.assert0}
\parbox{4in}{The group $G(k)$ acts transitively on the set of $d$-dimensional subspaces of $V$ satisfying the properties $\cP$.}
\end{equation}
Since the $\delta$-spaces are one $G(k)$-orbit by definition, this proves:
\begin{equation} \label{trans.assert1}
\parbox{4in}{The $\delta$-spaces are precisely the $d$-dimensional subspaces of $V$ satisfying the properties $\cP$.}
\end{equation}
 In many cases, we will refer to the literature for a proof of \eqref{trans.assert0}.  The proofs in the literature use various interpretations of the standard representation as the vector space underlying some algebraic structure.  For example, in the type $A$ example in \S\ref{A.eg}, we used the fact that $SL_n$ acts transitively on the subspaces of $k^n$ of a given dimension.
 
In fact, it suffices to prove \eqref{trans.assert0} in the case where $k$ is algebraically closed together with the statement
\begin{equation} \label{scalar}
\parbox{4in}{If $X$ satisfies $\cP$, then $X \ot \kbar$ satisfies $\cP$,}
\end{equation}
where $X$ is a $d$-dimensional subspace of $V$ and $\kbar$ is an algebraic closure of $k$.
To see this, suppose that $X$ satisfies $\cP$.  By our two hypotheses, there is some $g \in G(\kbar)$ such that $g (V_\delta \ot \kbar) = X \ot \kbar$.  For every $\s$ in the Galois group of $\kbar/k$, we find that 
\[
\s(g) (V_\delta \ot \kbar) = \s (X \ot \kbar) = X \ot \kbar = g (V_\delta \ot \kbar),
\]
hence the stabilizers  $\s(g) P_\delta \s(g)^{-1}$ and $gP_\delta g^{-1}$ agree.  That is, $gP_\delta g^{-1}$ is invariant under every element of in the Galois group, so it is defined over $k$ \cite[AG.14.4]{Borel} and there is some $h \in G(k)$ such that $hP_\delta h^{-1}$ equals $gP_\delta g^{-1}$ \cite[21.12]{Borel}.  This implies that $h (V_\delta \ot \kbar) = g (V_\delta \ot \kbar) = X \ot \kbar$, hence $h V_\delta$ equals $X$.

In the examples below, $\cP$ will almost always be a statement such as ``a certain polynomial vanishes'', for which property \eqref{scalar} clearly holds.  In such cases, we will omit any discussion of \eqref{scalar}.  The only exception will be for certain spaces related to the $E_6$ geometry, see \S\ref{E6}.

\medskip

In ``incidence'', we give a concrete description of how to tell if a $\delta$- and a $\delta'$-subspace are incident.  Our final description is purely in terms of subspaces of $V$, with no mention of the corresponding parabolic subgroups.  In most cases, a $\delta$- and $\delta'$-subspace will be incident if and only if one contains the other.  When this occurs, we will say that \emph{incidence is the same as inclusion}.
In this ``incidence'' portion, we return to the techniques of representation theory and eschew algebraic interpretations of the representation.  We do not need to consider the case where $\delta$ equals $\delta'$, because two spaces of the same type are incident if and only if they are equal.

\begin{rmk}  \label{min.free}
The hypothesis that $k$ has characteristic zero smoothes the presentation, but is in some cases unnecessary.  
(The rest of this remark will not be used elsewhere in the paper, so we omit many details.)  The material in Sections \ref{Tits} and \ref{realize} made no use of the characteristic zero hypothesis.  Regarding Section \ref{realize.2}, it remains true in all characteristics that there is a unique irreducible module $V$ with highest weight  $\omega$ for each dominant weight $\omega$, see e.g.\ \cite[24.4]{Borel}.  However, sometimes the weights and dimension of $V$ are not what one would expect coming from characteristic zero.  (For example, in characteristic 2 the standard representation of $B_n$ has dimension $2n$, not $2n + 1$.)  
On the positive side, Lemma~\ref{not.comp} remains unchanged.

Moreover, all the results of \S\ref{realize.2} hold \emph{regardless of the characteristic of $k$} when $V$ is the standard representation of $A_n$, $C_n$, $D_n$, $E_6$, or $E_7$.  In those cases, the highest weight $\omega$ is minuscule, i.e., all of the weights smaller than $\omega$ are Weyl-conjugates of $\omega$, hence the weights of $V$ are precisely the Weyl-orbit of $\omega$, regardless of the characteristic.  Moreover, the restricted weight $\wbar$ is a minuscule weight for $L_\delta$ for all $\delta \in \D \setminus \{ \beta \}$.
  \end{rmk}

\section{Example: type $D$ (orthogonal geometry)}  \label{Dn} \label{Dn.eg}

Consider the split simply connected group $G$ of type $D_n$ with $n \ge 4$.  This group is sometimes denoted $\Spin_{2n}$.  The geometry in this case is more complicated than for type $A$, apparently because the Dynkin diagram has a fork in it.  This case will illustrate the basic principles involved in handling forking diagrams, and we will use them when treating the $E$-groups later.  This geometry will be further investigated in Sections \ref{Dn.dual} and \ref{D4} below.

\begin{borel}{Dimensions and properties} \label{dim.D}
As in the type $A$ case, we compute the dimension of the $V_\delta$'s using Prop.~\ref{fund}.
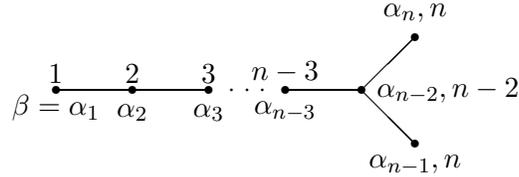
\[
\begin{picture}(7,2)
    \multiput(1,1)(1,0){2}{\line(1,0){1}}
    \put(4,1){\line(1,0){1}}
    
    \put(1,0.6){\makebox(0,0.4)[b]{$\beta = \alpha_1$}}  
    \put(2,0.6){\makebox(0,0.4)[b]{$\alpha_2$}}
    \put(3,0.6){\makebox(0,0.4)[b]{$\alpha_3$}}
    \put(4,0.6){\makebox(0,0.4)[b]{$\alpha_{n-3}$}}
    \put(1,1.1){\makebox(0,0.4)[b]{$1$}}
    \put(2,1.1){\makebox(0,0.4)[b]{$2$}}
    \put(3,1.1){\makebox(0,0.4)[b]{$3$}}
    \put(4,1.1){\makebox(0,0.4)[b]{${n-3}$}}
    \put(5.2,0.8){\makebox(0,0.4)[l]{$\alpha_{n-2}, n-2$}}

    \multiput(3.3,1)(0.2,0){3}{\circle*{0.01}}

    \put(5,1){\line(1,-1){0.7}}
    \put(5,1){\line(1,1){0.7}}

    \multiput(1,1)(1,0){5}{\circle*{\darkradG}}

    \put(5.7,0.3){\circle*{\darkradG}}
    \put(5.7,1.7){\circle*{\darkradG}}
    
    \put(5.7,-0.1){\makebox(0,0.4)[b]{$\alpha_{n-1}, n$}}
    \put(5.7,1.9){\makebox(0,0.4)[b]{$\alpha_n, n$}}
    
\end{picture}
\]

The standard representation of $G$ has dimension $2n$.  Since $-w_0 \omega_1$ is $\omega_1$, where $w_0$ is the longest element of the Weyl group, there is a nondegenerate $G$-invariant bilinear form $b$ on $V$, unique up to multiplication by an element of $k^\times$ \cite[8.7.5, Prop.~12]{Bou:g4}.  Moreover, $b$ is symmetric \cite[chap.~8, Table 1]{Bou:g4}.
For $\hw$ the highest weight vector in $V$ and $t$ an element of  the maximal torus $T$, we have:
   \[
   b(\hw,\hw) = b(t \cdot \hw, t\cdot \hw) = b(\omega(t) \hw, \omega(t) \hw) = \omega(t)^2 b(\hw, \hw).
   \]
Since $\omega$ is not the trivial character, $b(\hw, \hw)$ is 0.  Traditionally, a subspace $X$ is called \emph{isotropic} if $b(X,X)$ is zero.  We have just observed that the $\alpha_1$-spaces are isotropic.  By Prop.~\ref{incident.1}.1, the $\alpha_i$-spaces are isotropic for every $i$.
\end{borel}

\begin{borel}{Connection with the special orthogonal group} \label{SO}
In order to prove transitivity, we will now relate $G$ to the special orthogonal group $SO(b)$, the subgroup of $SL(V)$ preserving the bilinear form $b$.  Since $b$ is $G$-invariant and $G$ is connected, the representation of $G$ on $V$ is a homomorphism $\rho \!: G \ra SO(b)$; we claim that $\rho$ is a central isogeny.  Indeed, every proper, closed normal subgroup of $G$ is central, hence $\ker \rho$ is finite and the image of $\rho$ has the same dimension as $G$.  Root system data gives that the dimension of the Lie algebra of $G$ (equivalently, the dimension of $G$) is $\binom{2n}{2}$.  On the other hand, the Lie algebra of $SO(b)$ is isomorphic over an algebraic closure of $k$ to the space of skew-symmetric $2n$-by-$2n$ matrices, which also has dimension $\binom{2n}{2}$.  Since $\im \rho$ and $SO(b)$ are connected and have the same dimension, they are the same.  That is, $\rho$ is surjective.  The claim now follows because we are in characteristic zero, hence $\rho$ is automatically separable.
\end{borel}

\subsection*{Transitivity} 
We claim that $G$ acts transitively on the $m$-dimensional isotropic subspaces of $V$ for $m < n$.  (We refer the reader to \cite[Ch.~I]{Lam} for basic terminology and facts regarding symmetric bilinear forms.)  Let $X$, $X'$ be isotropic of dimension $m$.  They each lie in a direct sum of $m$ hyperbolic planes in $V$, and there is an isometry $f$ of $b$ that sends $X$ to $X'$ by Witt's Extension Theorem.  Since $V$ is isomorphic to a direct sum of $n$ hyperbolic planes, there is at least one plane where we may choose $f$ as we please.  If $f$ has determinant $-1$, we modify $f$ by a hyperplane reflection in this ``extra'' hyperbolic plane so that $f$ has determinant 1.  Over an algebraic closure $\kbar$ of $k$, $f$ is in the image of the map $G(\kbar) \ra SO(b)(\kbar)$, which proves the claim.  Moreover, we have proved that the $\alpha_i$-spaces are the $i$-dimensional isotropic subspaces for $1 \le i \le n - 2$.

We claim that the $n$-dimensional isotropic subspaces of $V$ make up two $G(k)$-orbits.  Since the stabilizers of $V_{\alpha_n}$ and $V_{\alpha_{n-1}}$ are the non-conjugate subgroups $P_{\alpha_n}$ and $P_{\alpha_{n-1}}$, there are at least two $G$-orbits.  On the other hand, consider the orthogonal group $O(b)$, i.e., the subgroup of $GL(V)$ preserving $b$; it is generated by $SO(b)$ and a hyperplane reflection.  In particular, $SO(b)$ is a normal subgroup of index 2.  An argument as in the preceding paragraph shows that $O(b)$ acts transitively on the set of $n$-dimensional isotropic subspaces, hence $SO(b)$ has at most two orbits.  We have shown that
  \[
  \{ \text{$\alpha_n$- and $\alpha_{n-1}$-spaces} \} = \{ \text{isotropic subspaces of dimension $n$} \}.
  \]
  
\begin{borel}{Incidence} \label{incidence.D}
Consider an $\alpha_i$-space $X_i$ and an $\alpha_j$-space $X_j$ with $i < j$.  If $(i, j)$ is not $(n - 1, n)$, then Prop.~\ref{incident.1}.2 applies and incidence is the same as inclusion.  We now argue that an $\alpha_{n-1}$-space and an $\alpha_n$-space are incident if and only if their intersection has dimension $n - 1$.

Let $M$ be the subgroup of $G$ generated by the $\X_\rho$ as $\rho$ ranges over the roots with $\alpha_{n-1}$- and $\alpha_n$-height zero; it has type $A_{n-2}$.  Write $Y$ for the subspace $M\hw$ of $V$; it is the standard representation of $M$ by the reasoning in \S\ref{realize.2} and so has dimension $n - 1$.  It is contained in $V_{\alpha_n}$, so it is isotropic.  

The parabolic $P_\anmn$ stabilizes $Y$ by the argument in the proof of Prop.~\ref{stable}.  The only subgroups of $G$ properly containing $P_\anmn$ are $P_{\alpha_{n-1}}$, $P_{\alpha_n}$, and $G$, none of which stabilize $Y$.  (For example, $P_{\alpha_n} Y$ contains $L_{\alpha_n} \hw$, which has dimension $n$.) Hence the stabilizer of $Y$ is precisely $P_\anmn$.

Suppose first that $X_{n-1}$ and $X_n$ are incident, so there is a $g \in G(k)$ such that $gX_{n-1} = V_{\alpha_{n-1}}$ and $gX_n = V_{\alpha_n}$.  The intersection $gX_{n-1} \cap gX_n$ contains $Y$ and so has dimension at least $n - 1$.  Because $X_{n-1}$ has dimension $n$ and is not contained in $X_n$, this shows that the intersection has dimension exactly $n - 1$.

Conversely, suppose that $X_{n-1} \cap X_n$ has dimension $n - 1$.  For $i = n - 1, n$, there is a $g_i \in G(k)$ such that $gX_i = V_{\alpha_i}$ and $g_i(X_{n-1} \cap X_n) = Y$ by the argument in the proof of Prop.~\ref{incident.1}(2).  Therefore the stabilizer $g_i P_{\alpha_i} g_i^{-1}$ of $X_i$ contains the stabilizer of $X_{n-1} \cap X_n$, which is a parabolic subgroup of $G$.  That is, $X_{n-1}$ and $X_n$ are incident.
\end{borel}

\begin{rmk}
A consequence of the above argument is that $V_{\alpha_{n-1}}$ and $V_{\alpha_n}$ are the unique $\alpha_{n-1}$ and $\alpha_n$ spaces containing $Y$.  Indeed, if $X$ is an $\alpha_i$-space containing $Y$ for $i = n - 1$ or $n$, then the stabilizer of $X$ contains $P_\anmn$ by the preceding paragraph, hence the stabilizer of $X$ is $P_{\alpha_i}$.  The bijection between objects in $\G_V$ and their stabilizers gives that $X$ equals $V_{\alpha_i}$.

Since $G(k)$ acts transitively on the $(n-1)$-dimensional isotropic subspaces, we have proved:
\emph{Each $(n-1)$-dimensional isotropic subspace is the intersection of two uniquely determined and incident $n$-dimensional subspaces, one of type $\alpha_n$ and one of type $\alpha_{n-1}$.}
This is a standard result in quadratic form theory, usually proved by quadratic-form-theoretic methods, see e.g.\ \cite[III.1.11]{Chev}.
\end{rmk}

\begin{borel}{An alternative view}
We now outline the geometry  that one obtains from $G$ by considering the fundamental representation with highest weight $\alpha_n$ (a ``half-spin'' representation) instead of the standard representation.  We continue with the same definitions of $V$, $V_{\alpha_i}$, etc., as in the rest of this section.  We may identify the vector space underlying the half-spin representation $S$ with $\wedge^{{\mathrm{even}}} V_{\alpha_n}$ as described in \cite[chap.~3]{Chev}.

For $i \ne n - 1$, the subspace $V_{\alpha_n}$ contains $V_{\alpha_i}$ and the ideal of $\wedge V_{\alpha_n}$ generated by $\wedge^i V_{\alpha_i}$ is stabilized by $P_{\alpha_i}$.  (Recall that the action of $G$ on $S$ is not precisely the standard action of $G$ on $\wedge V$, see \cite[\S2.2]{Chev}.  In particular, $G$ stabilizes $\wedge V_{\alpha_n}$.)  Following the naive algorithm in \S\ref{realize}, we take the intersection of $(\wedge^i V_{\alpha_i}) \wedge (\wedge V_{\alpha_n})$ with $S$ to be the subspace corresponding to $P_{\alpha_i}$.  Thinking in terms of exterior powers of vector spaces, it is clear that the $\alpha_i$-spaces in the half-spin representation have dimension $2^{n - i - 1}$ for $i \le n - 2$; they correspond to the  right ideals in the even Clifford algebra constructed in \cite[\S1]{G:flag}.  The $\alpha_n$-spaces are 1-dimensional and are the ``pure spinors'' corresponding to even maximal isotropic subspaces in the language of \cite[\S3.1]{Chev}.
We do not know how to describe the $\alpha_{n-1}$-spaces in this geometry.
\end{borel}

\section{Example: type $E_6$}  \label{E6}

The geometry for the split simply connected group $G$ of type $E_6$ exhibits two complexities.  We are prepared for the first---the fork in the diagram---thanks to our work in the previous section on groups of type $D$.  The second complication is new: the root $\delta = \alpha_6$ does not satisfy  the hypotheses of our workhorse Prop.~\ref{incident.1}.

This is the last example we will do in detail.  In the final section of this paper, we will give an explicit description of the duality in this geometry.

\subsection*{Dimensions and properties} The Dynkin diagram for $E_6$, labeled with the dimensions of the corresponding spaces is:
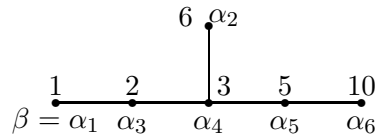
\[
\begin{picture}(7,2)
    \multiput(1,0.5)(1,0){5}{\circle*{\darkradG}}
    \put(3,1.5){\circle*{\darkradG}}

    \put(1,0.5){\line(1,0){4}}
    \put(3,1.5){\line(0,-1){1}}
    
    \put(1,0.1){\makebox(0,0.4)[b]{$\beta = \alpha_1$}}
    \put(2,0.1){\makebox(0,0.4)[b]{$\alpha_3$}}
    \put(3,0.1){\makebox(0,0.4)[b]{$\alpha_4$}}
    \put(4,0.1){\makebox(0,0.4)[b]{$\alpha_5$}}
    \put(5,0.1){\makebox(0,0.4)[b]{$\alpha_6$}}
    \put(3.2,1.5){\makebox(0,0.4)[b]{$\alpha_2$}}

    \put(1,.6){\makebox(0,0.4)[b]{$1$}}
    \put(2,.6){\makebox(0,0.4)[b]{$2$}}
    \put(3.2,.6){\makebox(0,0.4)[b]{$3$}}
    \put(4,.6){\makebox(0,0.4)[b]{$5$}}
    \put(5,.6){\makebox(0,0.4)[b]{$10$}}
    \put(2.7,1.5){\makebox(0,0.4)[b]{$6$}}
\end{picture}
\]

\begin{borel*} \label{min}
We will make a detailed study of the weights of representation $V$ in order to understand  the algebraic properties satisfied by the objects in the geometry.  Recall from, e.g., \cite[\S{VIII.7.3}]{Bou:g4} that a weight $\la$ is \emph{minuscule} if $\qform{w \la, \alpha} = 1$, 0, or $-1$ for every element $w$ of the Weyl group and every root $\alpha$.  Further, $\omega_1$ is a minuscule weight for $E_6$, and this implies that the weights of $V$ are those weights in the Weyl orbt of $\omega_1$ and all have multiplicity one.  

Figure \ref{weights} displays the 27 weights of $V$ as a Hasse diagram relative to the usual partial ordering of the weights (recalled in \ref{spindle}).  An edge joining two weights $\la > \mu$ is labeled with $i$ if $\la - \alpha_i = \mu$.  The row vectors list the coordinates of the weights with respect to the basis consisting of fundamental weights.    The lowest weight of $V_{\alpha_i}$ (cf.~\ref{spindle}) is labeled $\la_i$.  (This diagram is an elaboration of \cite[Fig.~20]{PSV}.  To construct such a diagram from scratch, one calculates the Weyl orbit of $\omega_1$---this amounts to a single instruction in the software packages LiE \cite{LiE} or Magma \cite{Magma}.  Alternatively, Dynkin \cite[p.~333]{Dynkin} gives a general algorithm for finding 
the weights of an irreducible representation from its highest weight using \cite[VIII.7.2, Prop.~3(i)]{Bou:g4}.
Next it is helpful to observe that all the weights are of the form $\omega_1 - \alpha$ for a sum of positive roots $\alpha$, so the weights are stratified by the height of $\alpha$.)

We remark that one can see from Figure \ref{weights} that $\omega_1$ is minuscule; indeed, every coordinate of every weight is 1, 0, or $-1$, as required.  We will make repeated use of this fact in \S\ref{E6.dual}.
\end{borel*}

\medskip

\renewcommand{\thefigure}{\theequation}
\refstepcounter{equation} \label{weights}
\begin{figure}[th]
\begin{minipage}{\textwidth}
\[
\xymatrix@R=8pt{
& {(1, 0, 0, 0, 0, 0) = \omega_1 = \la_1} \ar@{-}[d]^1 & \\
& (-1, 0, 1, 0, 0, 0) = \la_3 \ar@{-}[d]^3 & \\
& (0, 0, -1, 1, 0, 0) = \la_4 \ar@{-}[d]^4 & \\
& (0, 1, 0, -1, 1, 0) \ar@{-} [ld]_2 \ar@{-}[rd]^5 & \\
\la_5 = (0, -1, 0, 0, 1, 0) \ar@{-}[d]_5 & & (0, 1, 0, 0, -1, 1)  \ar@{-}[dll]^2 \ar@{-}[d]^6\\
(0, -1, 0, 1, -1, 1) \ar@{-}[d]_4 \ar@{-}[drr]_6 & & (0, 1, 0, 0, 0, -1) = \la_2 \ar@{-}[d]^2 \\
(0, 0, 1, -1, 0, 1) \ar@{-}[d]_3 \ar@{-}[drr]_6 & & (0, -1, 0, 1, 0, -1) \ar@{-}[d]^4 \\
(1, 0, -1, 0, 0, 1) \ar@{-}[d]_1 \ar@{-}[dr]^6 && (0, 0, 1, -1, 1, -1) \ar@{-}[d]^5 \ar@{-}[dl]_3 \\
\la_6 = (-1, 0, 0, 0, 0, 1)  \ar@{-}[d]_6 & (1, 0, -1, 0, 1, -1)\ar@{-}[dr]_5 \ar@{-}[dl]^1 & (0, 0, 1, 0, -1, 0) \ar@{-}[d]^3 \\
(-1, 0, 0, 0, 1, -1) \ar@{-}[d]_5 & & (1, 0, -1, 1, -1, 0) \ar@{-}[d]^4 \ar@{-}[dll]_1\\
(-1, 0, 0, 1, -1, 0) \ar@{-}[d]_4& & (1, 1, 0, -1, 0, 0)\ar@{-}[d]^2 \ar@{-}[dll]_1 \\
(-1, 1, 1, -1, 0, 0) \ar@{-}[d]_3 \ar@{-}[drr]^2 & & (1, -1, 0, 0, 0, 0)\ar@{-}[d]^1 \\
(0, 1, -1, 0, 0, 0) \ar@{-}[dr]_2&& (-1, 1, 1, 0, 0, 0) \ar@{-}[dl]^3 \\
& (0, -1, -1, 1, 0, 0) \ar@{-}[d]^4& \\
& (0, 0, 0, -1, 1, 0) \ar@{-}[d]^5& \\
& (0, 0, 0, 0, -1, 1) \ar@{-}[d]^6& \\
& (0, 0, 0, 0, 0, -1) = -\omega_6 
}
\]
\end{minipage}
\caption{Hasse diagram of the weights of $V$, where $(c_1, c_2, \ldots, c_6)$ denotes the weight $\sum_{i=1}^6 c_i \omega_i$.  The map $-\phi$ reflects the diagram across its horizontal axis of symmetry.}
\end{figure}
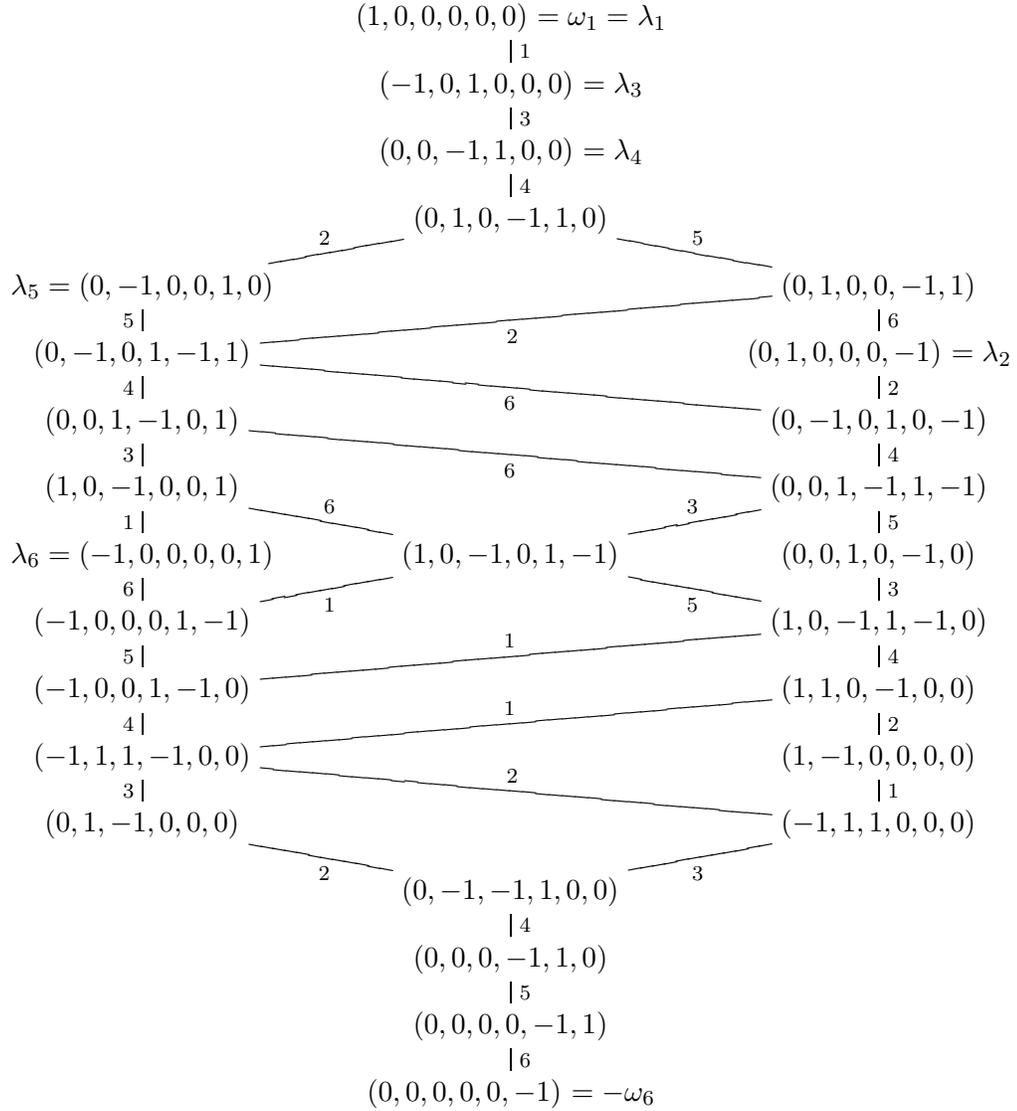

The automorphism of order 2 of the Dynkin diagram gives an automorphism of the root system, hence an automorphism $\phi$ of $G$ of order 2.  

\begin{lem}
The subgroup of $G$ consisting of elements fixed by $\phi$ is split simple of type $F_4$.
\end{lem}

The lemma is well-known.  We  sketch a proof for the reader's convenience.

\begin{proof}[Sketch of proof]
Steinberg \cite{St:end} gives that the fixed subgroup $G_\phi$ is connected (his 9.7) and reductive (his 9.4).  We now inspect its Lie algebra $\g_\phi$, which consists of the elements of the Lie algebra $\g$ of $G$ that are fixed by $\phi$.  Fix a Chevalley basis $\{ x_\alpha \mid \alpha \in \Phi \} \cup \{ h_\delta \mid \delta \in \D \}$ for $\mathfrak{g}$ as in \cite[\S8]{Hum:LA}.  The elements $x_\alpha + \phi(x_\alpha)$ for $\alpha \in \Phi$ and $h_\delta + \phi(h_\delta)$  for $\delta \in \D$ span $\g_\phi$ since $\phi$ has order 2.  Moreover, the $h_\delta$'s are uniquely determined, hence $\phi(h_\alpha)$ equals $h_{\phi(\alpha)}$ for all $\alpha \in \Phi$.  Computing the root space decomposition of $\g_\phi$ with respect to the torus spanned by
\begin{equation} \label{F4.tor}
h_{\alpha_1} + h_{\alpha_6}, \ h_{\alpha_3} + h_{\alpha_5}, \ h_{\alpha_4}\text{, and } h_{\alpha_2},
\end{equation}
we find that $\g_\phi$ is simple of type $F_4$, where the displayed elements correspond to the coroots $\ach_4$, $\ach_3$, $\ach_2$, and $\ach_1$ respectively.
\end{proof}

For simplicity, we denote the fixed subgroup by $F_4$.  We can compute the restriction of weights of $V$ to $F_4$ on the level of Lie algebras:  A weight $\la := \sum c_i \omega_i$ maps $h_{\alpha_j} \mapsto c_j$, hence it maps the elements listed in \eqref{F4.tor} to 
\[
(c_1 + c_6),\ (c_3 + c_5),\ c_4\text{, and } c_2
\]
respectively.
In terms of the fundamental weights of $F_4$, the restriction of $\la$ is $c_2 \omega_1 + c_4 \omega_2 + (c_3+ c_5)\omega_3 + (c_1 + c_6) \omega_4$.
We find that $V$ decomposes (as a representation of $F_4$) as a direct sum of a 1-dimensional trivial representation $C$ and the standard representation of $F_4$, which we denote by $V_0$.

\begin{prop} \label{bil.prop}
There is a bilinear form $b$ on $V$ such that
\begin{equation} \label{b.prop}
b(\phi(g)x, gy) = b(x, y) \quad \text{for all $g \in G$ and $x, y \in V$.}
\end{equation}
It is unique up to multiplication by a scalar.  Moreover, it is symmetric and nondegenerate, and $b|_C$ is not zero.
\end{prop}

\begin{proof}
First we construct a bilinear form $b$ on $V$ satisfying \eqref{b.prop}.  Write $\rho \!: G \ra GL(V)$ for the representation of $G$ on $V$, and write $\rho^* \!: G \ra GL(V^*)$ for the dual representation defined by
  \[
  (\rho^*(g) f) (x) := f(\rho(g)^{-1} x) \quad \text{for $g \in G$, $f \in V^*$, and $x \in V$.}
  \]
The representations $\rho \phi$ and $\rho^*$ are both irreducible with highest weight $\omega_6$, hence they are isomorphic.  Fix an isomorphism $h \!: V \ra V^*$ such that 
$h \rho \phi(g) h^{-1} = \rho^*(g)$ for all $g \in G$.  Define $b$ by setting
   \[
   b(x, y) := h(x)(y).
   \]
This $b$ is clearly bilinear and 
  \[
  b(\phi(g) x, g y) = h(\rho\phi(g)x) (gy) = \left[ \rho^*(g) h(x) \right] (gy) = b(x, y).
  \]

We now argue that any bilinear form $b$ satisfying \eqref{b.prop} is symmetric.  Set 
\[
b_\e(x, y) := b(x, y) + \e b(y, x).
\]
Then $b_1$ and $b_{-1}$ are bilinear, $b_1$ is symmetric, 
$b_{-1}$ is skew-symmetric, and $2b = b_1 + b_{-1}$.  We prove that $b_{-1}$ is identically zero.  In any case, $b_{-1}$ satisfies \eqref{b.prop} (using that $\phi$ is its own inverse), hence $b_{-1}$ is $F_4$-invariant.  But $V_0$ does not support a nonzero $F_4$-invariant skew-symmetric form,
hence $b_{-1}$ restricts to zero on $V_0$.  Fix $x \in V_0$ a nonzero vector with a nonzero weight $\la$ with respect to $F_4$, and let $c$ be a nonzero vector in $C$.  Since
   \[
   b_\e(c, x) = b_\e(tc, tx) = \la(t) b_\e(c, x) 
   \]
for $t$ in the $F_4$-torus, $b_\e(c,x)$ is zero.  Since $V_0$ is an irreducible representation of $F_4$, $b_\e(c, V_0)$ is zero.  But $b_{-1}$ is skew-symmetric, so $b_{-1}(c, c)$ is also zero, and we have proved the claim.

The previous paragraph also gives more.  Continue the assumption that $b$ satisfies \eqref{b.prop} and suppose that $b$ is not identically zero.
Then $C$ and $V_0$ are orthogonal subspaces.  For $x \in V$ and $r$ in the radical of $b$, we have
   \[
   b(gr, x) = b(r, \phi(g)^{-1}x) = 0,
   \]
hence the radical is $G$-invariant.  Since $V$ is irreducible and $b$ is not identically zero, the radical is zero, i.e., $b$ is nondegenerate.
Since $C$ is 1-dimensional and orthogonal to $V_0$, we find that $b$ restricts to be nonzero on $C$.

We now prove uniqueness.  Let $b$, $b'$ be bilinear forms on $V$ satisfying \eqref{b.prop}.  The representation $V_0$ of $F_4$ supports a unique symmetric bilinear form up to a scalar multiple, so by modifying $b'$ by a factor in $k^\times$, we may assume that $b$ and $b'$ have the same restriction to $V_0$.  Then $b - b'$ is a bilinear form on $V$ satisfying \eqref{b.prop} that restricts to be zero on $V_0$.  By the previous paragraph, $b - b'$ is identically zero, and we have proved uniqueness.
\end{proof}

\begin{borel*} \label{repth}
We can use representation theory to find $G$-invariant polynomial functions on $V$.  Plugging the formal character of the dual representation $V^*$ into the degree $d$ complete symmetric polynomial in $\dim V^*$ variables gives the formal character of $S^d(V^*)$, the $d$-th symmetric power of $V^*$.  From this, one can write $S^d(V^*)$ as a direct sum of irreducible representations \cite[22.5A]{Hum:LA}.  (For small $d$, 
these computations are easily done using a computer package like LiE or Magma.)  If $S^d(V^*)$ has a unique 1-dimensional summand---corresponding to a summand with highest weight $0$---then $V$ supports a $G$-invariant homogeneous polynomial of degree $d$, uniquely determined up to a factor in $k^\times$.  This happens for $G$ of type $E_6$ and $d = 3$.

Write $N$ for a nonzero cubic form on $V$ as discovered in the previous paragraph.
We abuse notation and write $N$ also for the trilinearization of $N$ on $V$ such that
$N(x, x, x) = 6N(x)$ for all $x \in V$. 
Let $\#$ denote the bilinear product defined implicitly by the formula
  \begin{equation} \label{sharp.def}
  b(x \# y, z) = N(x, y, z) \quad \text{for $x, y, z \in V$.}
  \end{equation}
Using \eqref{b.prop}, we find
  \begin{equation} \label{duality}
  \phi(g) (x \# y) = (gx) \# (gy) \quad \text{for $g \in G$ and $x, y \in V$.}
  \end{equation}
We write $x^\#$ for $(x \# x)/2$.  (The factor 6 above arises naturally from multilinearization, see \cite[\S{IV.5.8}, Prop.~12(i)]{Bou:alg2}.  The factor 2 here occurs for the same reason:  In order to apply identities from Jordan theory in \S\ref{E6.dual}, we adopt the Jordan theorist view that the square $x \mapsto x^\#$ is the base object and the bilinear product is obtained by multilinearizing it.)
\end{borel*}

\medskip

\begin{borel*} \label{singular}
The same argument as in \ref{dim.D} shows that the $\alpha_1$-spaces are 1-dimensional subspaces consisting of elements $x \in V$ such that $x^\# = 0$.  We say a nonzero vector $x \in V$ is \emph{singular} if $x^\# = 0$.  (These 1-dimensional subspaces are precisely the singular points for the hypersurface in $\mathbb{P}(V)$ defined by $N = 0$, because $b(x^\#, y)$ is the directional derivative of $N$ at $x$ in the direction $y$.)  We call a subspace of $V$ singular if its nonzero elements are singular.  By Prop.~\ref{incident.1}.1, the $\alpha_i$-spaces are singular for $i \ne 6$.
\end{borel*}

We will now  investigate the restriction of the representation of $G$ on $V$ to the subgroup $L_{\alpha_6}$ of type $D_5$.  This will give us finer information about the product $\#$ and lead us to a description of the $\alpha_6$-spaces.  
To see how a weight of $G$ restricts to $L_{\alpha_6}$, one drops the last coordinate and moves the second coordinate to the end of the vector (to allow for the fact that weights of $D_5$ and $E_6$ are numbered somewhat incompatibly in \cite{Bou:g4}).  

Let $W'$ be the subgroup of the Weyl group $W$ of $G$ generated by the reflections with respect to the roots $\alpha_i$ for $i \ne 6$.  It is the Weyl group of $L_{\alpha_6}$, and it is the stabilizer of $-\omega_6$ in $W$ \cite[Th.~1.12c]{Hum:ref}.

\begin{lem} \label{orbits}
The orbits of $W'$ in the weights of $V$ are 
the weights $\ge \la_6$, the weights between $\la_2$ and $(0, 0, 0, 0, -1, 1)$, and the weight $-\omega_6$.
\end{lem}

\begin{proof}
Since the highest weight $\omega_1$ of $V$ is minuscule, we have $\qform{\mu, \alpha} = 1$, $0$ or $-1$ for every weight $\mu$ of $V$ and every root $\alpha$.  If $\mu$ and $\mu - \delta$ are both weights for some $\delta \in \D$, then $\qform{\mu, \delta} = 1$, $\qform{\mu - \delta, \delta} = -1$, and the reflection $s_\delta$ with respect to the root $\delta$ interchanges $\mu$ and $\mu - \delta$.  Consulting Figure \ref{weights}, we see that $W'$ acts transitively on each of the three sets of weights named in the statement of the lemma.

Conversely, $\omega_1$ and $\la_2$ restrict to the weights $(1, 0, 0, 0, 0)$ and $(0, 0, 0, 0, 1)$ on $L_{\alpha_6}$, which are not congruent modulo the $D_5$ root lattice.  Therefore, they lie in different $W'$-orbits \cite[VI.1.9, Prop.~27]{Bou:g4}.
\end{proof}

By restricting the weights of $V$ to $L_{\alpha_6}$, we can decompose $V$ as a direct sum of irreducible representations.  The proof of Lemma \ref{orbits} shows that the components of $V$ are
\begin{itemize}
\item the standard representation $V_{\alpha_6}$ of $L_{\alpha_6}$ (with highest weight $(1, 0, 0, 0,0)$),
\item a half-spin representation (with highest weight $(0, 0, 0, 0, 1)$), and
\item a 1-dimensional trivial representation (from the lowest weight vector $-\omega_6$).
\end{itemize}

\begin{cor} \label{triples}
The Weyl group of type $E_6$ acts transitively on triples $\mu_1, \mu_2, \mu_3$ of weights of $V$ such that $\mu_1 + \mu_2 + \mu_3 = 0$.
\end{cor}

\begin{proof}
The Weyl group acts transitively on the weights of $V$, so we may assume that $\mu_1$ is $-\omega_6$.  Since $\mu_3 = \omega_6 - \mu_2$ is a weight, $\mu_2$ cannot have last coordinate equal to $-1$, otherwise $\mu_3$ would have last coordinate $-2$, which is impossible.  In particular, $\mu_2$ cannot be $-\omega_6$ or $\la_2$.  Since the set of triples $-\omega_6, \mu_2, \mu_3$ with sum 0 is stable under the action of $W'$, $\mu_2$ must lie in the $W'$-orbit with lowest weight $\la_6$.
\end{proof}

Our preliminary results about the action of the Weyl group can now give us concrete information about the product $\#$.

\begin{lem} \label{nonzero}
Let $x_1, x_2$ be nonzero vectors in $V$ of weight $\mu_1, \mu_2$ respectively.  The product $x_1 \# x_2$ is nonzero if and only if $\phi(\mu_1 + \mu_2)$ is a weight of $V$.
\end{lem}

The equations
\begin{equation} \label{lambdas}
\phi(\la_2 + \la_5) = \la_3 \quad \text{and} \quad \phi(\omega_1 + \la_6) = \omega_1
\end{equation}
furnish specific examples where the product $\#$ is not zero.

\begin{proof}[Proof of Lemma \ref{nonzero}]
If $\phi(\mu_1 + \mu_2)$ is not a weight of $V$, then the product is zero by \eqref{duality}.  So suppose that $\phi(\mu_1 + \mu_2)$ is a weight of $V$.  Since $-\phi$ is in the Weyl group, $\mu_3 := -\mu_1 - \mu_2$ is a weight of $V$; let $y$ be a nonzero vector with that weight.  

We claim that $N(x_1, x_2, y)$ is not zero, and hence that $x_1 \# x_2$ is not zero.  Recall that every weight of $V$ has multiplicity 1 because $V$ is minuscule.  Fix a basis $\{ b_\la \}$ for $V$ where $b_\la$ has weight $\la$, and write $N$ in terms of the dual basis $\{ f_\la \}$.  Since $N$ is $G$-invariant, a monomial $f_{\nu_1} f_{\nu_2} f_{\nu_3}$ has zero coefficient if $\nu_1 + \nu_2 + \nu_3$ is not zero.  On the other hand, since $N$ is not identically zero, there exist weights $\nu_1, \nu_2, \nu_3$ such that the coefficient of $f_{\nu_1} f_{\nu_2} f_{\nu_3}$ is not zero.  Since their sum $\nu_1 + \nu_2 + \nu_3$ is zero, there is an element $w$ in the Weyl group such that $w \nu_i = \mu_i$ for each $i$ by Cor.~\ref{triples}. A representative of $w$ can be found in $G$, hence the coefficient of $f_{\mu_1} f_{\mu_2} f_{\mu_3}$ in $N$ is not zero.  In particular, $N(x_1, x_2, y)$ is not zero, as claimed.
\end{proof}

Finally, we can give an explicit description of the $\alpha_6$-space $V_{\alpha_6}$.

\begin{prop} \label{hyper.lem}
$V_{\alpha_6} = \hw \# V$, where $\hw$ is the highest weight vector.
\end{prop}

\begin{proof}
($\supseteq$): Suppose that $x \in V$ has weight $\mu$, which is necessarily at least the minimum weight $-\omega_6$.   Since $\phi$ respects the partial ordering on the weights,  $\hw \# x$ has weight at least $\phi(\omega_1 - \omega_6) = \la_6$.  But this is the lowest weight of $V_{\alpha_6}$, so $V_{\alpha_6}$ contains $\hw \# x$ by \ref{spindle}.  Since every element of $V$ is a sum of weight vectors, we have shown that $V_{\alpha_6}$ contains $\hw \# V$.

($\subseteq$): Let $v \in V_{\alpha_6}$ be nonzero with weight $\mu$.  By Lemma \ref{orbits}, $\mu$ equals $w \la_6$ for some $w \in W'$; fix an element $n \in L_{\alpha_6}$ representing $w$.  
For $\lw$ a vector of the lowest weight $-\omega_6$, the vector $n(\hw \# \lw)$ is nonzero (by Lemma \ref{nonzero}) of weight $\mu$.  Since all weight spaces in $V$ are one-dimensional, $n(\hw \# \lw)$ is a nonzero multiple of $v$.  This proves that $V_{\alpha_6}$ is contained in $\hw \# V$.
\end{proof}

For $x$ singular in  $V$, we call the subspace $x \# V$ a \emph{hyperline}, following the terminology from \cite[p.~25]{Ti:Rsp} and \cite[p.~606]{Scharlau}.
Combining the proposition with \eqref{duality}, we find that every $\alpha_6$-space is a hyperline.

\begin{rmk}[Cf.~\protect{\cite[p.~17]{Faulk:oct}}] \label{hyper.quad} 
As the standard representation of the group $L_{\alpha_6}$ of type $D_5$, the space $V_{\alpha_6}$ supports an $L_{\alpha_6}$-invariant quadratic form, uniquely determined up to a scalar.   We claim that it is the form $q$ given implicitly by the equation
   \begin{equation} \label{qf.recipe}
   x^\# = q(x) \hw \quad \text{for $x \in V_{\alpha_6}$.}
   \end{equation}
   
First, observe that for $y, z$ weight vectors in $V_{\alpha_6}$, $y \# z$ has weight at least   
\[   
   \phi(2 \la_6) = \omega_1 - (\alpha_2 + \alpha_3 + 2 \alpha_4 + 2\alpha_5 + 2\alpha_6).
\] 
But the only edge leaving $\omega_1$ in Fig.~\ref{weights} is $\alpha_1$, so $\omega_1$ is the only possible weight for $y \# z$.
Therefore, for every $x \in V_{\alpha_6}$, the vector $x^\#$ is in the $k$-span of $\hw$ and the recipe \eqref{qf.recipe} defines a quadratic form on $V_{\alpha_6}$.

For $g \in L_{\alpha_6}$, we have
$q(gx) \hw = q(x) \phi(g) \hw$.  Since $\phi(g)$ fixes $\hw$ (by \ref{not.comp}, if you like), the form $q$ is $L_{\alpha_6}$-invariant.  Moreover, $q$ is not the zero form since $\#$ is bilinear and $\hw \# x$ is not zero for $x$ of weight $\la_6$, by Lemma \ref{nonzero}. 
Combining the two previous sentences, $q$ is the unique invariant quadratic form as claimed, and it is nondegenerate.  One consequence of this is that  $\hw \# V$ cannot contain an isotropic 6-dimensional subspace.  
In terms of the geometry, an $\alpha_6$-space cannot contain an $\alpha_2$-space.  A second consequence is that $(\hw \# V)^\#$ is $k\hw$; it follows that every $\alpha_6$-space $X$ equals $X^\# \# V$.
\end{rmk}

\begin{borel}{Connections with Jordan algebras} \label{conjord}
It is well-known that the cubic form $N$ can be realized as the generic norm (``determinant") on the split Albert algebra $A$, see e.g.\ \cite[7.3.1]{Sp:ex}.  (Recall that an Albert algebra is a 27-dimensional exceptional Jordan algebra, see e.g.\ \cite[Ch.~5]{Sp:ex}.)  That is, there is a vector space isomorphism $f\!: V \ra A$ such that the generic norm $N_A$ on $A$ satisfies 
\begin{equation} \label{N.eq}
N_A(f(x)) = N(x)  \quad \text{for $x \in V$.}
\end{equation}
We now show that we can pick $f$ such that \eqref{N.eq} holds and such that $f$ identifies $b$ with the trace bilinear form $b_A$ on $A$.  This turns out to be a bit technical, but we need it in order to apply results about Albert algebras in the ``transitivity'' portion.  We assume that $k$ is algebraically closed.

Write $b'$ for the symmetric bilinear form on $V$ induced by $b_A$ via $f$.  Formula \eqref{b.prop} with $b'$ instead of $b$ defines an automorphism $\phi'$ of $GL(V)$.  This $\phi'$ restricts to an automorphism of $G$ of order 2 by \cite[p.~76]{Jac:J3}, and it leaves a Borel subgroup $B'$ and a maximal torus $T'$ invariant by \cite[7.5]{St:end}.  There is some $h \in G(k)$ such that $hTh^{-1} = T'$ and $hBh^{-1} = B'$.  Replacing $f$ with the map $x \mapsto f(hx)$ replaces $\phi'$ with the map
   \[
   g \mapsto h^{-1} \phi'(hgh^{-1}) h.
   \]
Changing $f$ in this way, we may assume that $\phi$ and $\phi'$ both leave $T$ and $B$ invariant.  

Equation \eqref{b.prop} implies that $\phi$ and $\phi'$ have the same action on the center of $G$, which consists of cube roots of unity. Therefore, $\phi'(g) = t \phi(g) t^{-1}$ for some $t \in T(k)$ \cite[14.9]{Borel}.  Elementary computations as in \cite[2.7]{KMRT} show that $\phi(t) = t^{-1}$.  Since $k$ is algebraically closed, there is a ``square root'' $s$ of $t$, i.e., an element $s \in T(k)$ such that $s^2 = t$ and $\phi(s) = s^{-1}$.  Replacing $f$ with the map $x \mapsto f(sx)$, we may assume that $\phi'$ equals $\phi$.  Prop.~\ref{bil.prop} gives that $b'$ is a scalar multiple of $b$.  But $b$ was only defined up to a scalar multiple to begin with, so we may replace $b$ with $b'$ and we have found an isomorphism $f$ as desired.
\end{borel}   



\subsection*{Transitivity}
It is known that $G$ acts transitively on the $i$-dimensional singular subspaces for $i = 1, 2, 3, 4$, and $6$ by \cite[3.12]{SpV},  \cite[p.~33]{Faulk:oct}, \cite[p.~35, Cor.~5]{Racine}, or \cite[6.5(2)]{Asch:E6}.  Thus the $\alpha_1$-, $\alpha_2$-, $\alpha_3$-, and $\alpha_4$-spaces are as described in Table \ref{E6.table} below.
Every hyperline is by definition of the form $x \# V$ for a singular $x \in V$.  Since $G$ acts transitively on the 1-dimensional singular subspaces, it acts transitively on the hyperlines by \eqref{duality}.  Therefore, the $\alpha_6$-spaces are the hyperlines.

\smallskip

We now prove that the $\alpha_5$-spaces are the 5-dimensional, maximal  singular subspaces.  In the notation of \eqref{trans.assert0}--\eqref{scalar}, we take $d = 5$ and $\cP$ to be ``maximal singular". 
We have to be a little careful; for example, it is not clear that property \eqref{scalar} holds, that is, that the property ``maximal singular subspace'' is preserved when one goes from $k$ to $\kbar$.

As in \ref{incidence.D}, take $M$ to be the subgroup of $G$ generated by $\X_\rho$ as $\rho$ varies over the roots of $\alpha_2$- and $\alpha_6$-height zero; it is of type $A_4$.  Put $Y := M \hw$; it is the standard representation of $M$ and has dimension 5.  It is contained in $V_{\alpha_2}$, so it is singular.  The arguments in \ref{incidence.D} show that the stabilizer of $Y$ is $P_{ \{ \alpha_2, \alpha_6 \}}$.

Suppose first that $k$ is algebraically closed.
By \cite[p.~35, Cor.~5]{Racine}, the collection of 5-dimensional singular subspaces of $V$ is a union of two $G$-orbits, corresponding to those singular subspaces that are and are not maximal.  Since $V_5$ and $Y$ have non-conjugate stabilizers $P_{\alpha_5}$ and $P_{ \{ \alpha_2, \alpha_6 \}}$, they lie in different orbits.  This proves  \eqref{trans.assert0} for $k$ algebraically closed.

Now we treat the case where $k$ is arbitrary by proving \eqref{scalar} for the $\alpha_5$-spaces.  That is, let $X$ be a 5-dimensional maximal singular subspace of $V$.  For sake of contradiction, suppose that $X \ot \kbar$ is contained in a 6-dimensional singular subspace $\Zbar$ of $V \ot \kbar$.  We note that the arguments in \ref{incidence.D} show that $V_{\alpha_2} \ot \kbar$ is the unique $\alpha_2$-space containing $Y \ot \kbar$; since $G(\kbar)$ acts transitively on the 5-dimensional nonmaximal singular subspaces, every such subspace (e.g., $X \ot \kbar$) is contained in a unique 6-dimensional singular subspace.  For every $\s$ in the Galois group of $\kbar/k$, we have $X \ot \kbar = \s(X \ot \kbar) \subseteq \s(\Zbar)$. The uniqueness of $\Zbar$ implies that $\Zbar$ is stable under the Galois group, hence $\Zbar$ equals $Z \ot \kbar$ for some singular subspace $Z$ of $V$ \cite[AG.14.2]{Borel}.  That is, $X$ is not a maximal singular subspace of $V$.  This proves \eqref{scalar} for $X$ of type  $\alpha_5$, hence also the claimed description of the $\alpha_5$-spaces.

We summarize the descriptions of the $\alpha_i$-spaces in Table \ref{E6.table}.

\renewcommand{\thetable}{\theequation}
\refstepcounter{equation} \label{E6.table}
\begin{table}[th]
\begin{minipage}{\textwidth}
\begin{center}
\begin{tabular}{c|p{1.9in}p{1.6in}}
$\underline{\quad}$-space &description&name in \cite{SpV} \\ \hline
$\alpha_1$&$1$-dim'l singular& point\\
$\alpha_2$&$6$-dim'l singular& max'l space of 2nd kind\\
$\alpha_3$&$2$-dim'l singular& space of prdim 1\\
$\alpha_4$&$3$-dim'l singular& space of prdim 2\\
$\alpha_5$&$5$-dim'l, maximal singular& max'l space of 1st kind\\
$\alpha_6$&hyperline&line
\end{tabular}
\medskip
\caption{$\alpha_i$-spaces in the $E_6$ geometry}
\end{center}
\end{minipage}
\end{table}

\begin{borel}{Incidence} \label{incidence.E}
Let $X'$ be an $\alpha_i$-space and let $X$ be an $\alpha_j$-space with $i < j$.  If $(i, j)$ is not $(2, 5)$ or $(2, 6)$, Prop.~\ref{incident.1}.2 applies and incidence is the same as inclusion.
As in the $D_n$ case, we quickly find: An $\alpha_2$- and an $\alpha_5$-space are incident if and only if their intersection is 4-dimensional (equivalently, the stabilizer of their intersection is a parabolic of type $\{ \alpha_2, \alpha_5 \}$).  An $\alpha_2$- and an $\alpha_6$-space are incident if and only if their intersection is 5-dimensional (equivalently, the stabilizer of their intersection is a parabolic of type $\{ \alpha_2, \alpha_6 \}$).
\end{borel}

\begin{brmk}
We deduced the existence of a $G$-invariant cubic form on $V$ purely by computations with characters, but the invariant cubic form was written down long before these mathematical tools were available.  For example, it is in Cartan's 1894 thesis, see \cite[p.~143]{Car:th}.  Dickson pointed out in \cite{Dickson:bull} that the variables occurring in the cubic form correspond to the the 27 lines on a cubic surface.  For a modern explanation of the correspondence between the 27 lines on a cubic surface and  the 27 weights in Figure \ref{weights}, see \cite[\S\S23, 25]{Manin} or \cite[pp.~235, 236]{IS}.  Neher \cite[\S{II.5}]{Neher} describes the correspondence from a Jordan-theoretic viewpoint.  Lurie gives a derivation of the explicit formula for the cubic form in \cite{Lurie}.

Diagrams like Figure \ref{weights} for different groups and representations can be found in, for example, \cite{PSV}, \cite[\S2.5]{Scharlau}, and the references listed on p.~519 of \cite{Scharlau}.  
\end{brmk}

\section{Example: type $E_7$}  \label{E7}

We will now treat the split group of type $E_7$.  In the interest of brevity, this section is less detailed than the preceding ones.

\subsection*{Dimensions and properties} For $G$ of type $E_7$, the Dynkin diagram looks like 
\begin{equation} \label{E7.dynk}
\begin{picture}(7,1.7)
    \multiput(1,.5)(1,0){6}{\circle*{\darkradG}}
    \put(4,1.5){\circle*{\darkradG}}

    \put(1,.5){\line(1,0){5}}
    \put(4,1.5){\line(0,-1){1}}
    
    \put(1,0.1){\makebox(0,0.4)[b]{$\beta = \alpha_7$}}
    \put(2,0.1){\makebox(0,0.4)[b]{$\alpha_6$}}
    \put(3,0.1){\makebox(0,0.4)[b]{$\alpha_5$}}
    \put(4,0.1){\makebox(0,0.4)[b]{$\alpha_4$}}
    \put(5,0.1){\makebox(0,0.4)[b]{$\alpha_3$}}
    \put(6,0.1){\makebox(0,0.4)[b]{$\alpha_1$}}
    \put(4.4,1.4){\makebox(0,0.4)[b]{$\alpha_2$}}

    \put(1,.6){\makebox(0,0.4)[b]{$1$}}
    \put(2,.6){\makebox(0,0.4)[b]{$2$}}
    \put(3,.6){\makebox(0,0.4)[b]{$3$}}
    \put(3.8,.6){\makebox(0,0.4)[b]{$4$}}
    \put(5,.6){\makebox(0,0.4)[b]{$6$}}
    \put(6,.6){\makebox(0,0.4)[b]{$12$}}
    \put(3.8,1.4){\makebox(0,0.4)[b]{$7$}}

\end{picture}
\end{equation}
As in \ref{repth}, we find that there is a $G$-invariant quartic form $q$ and a skew-symmetric bilinear form $b$ on $V$; both are unique up to a factor in $k^\times$.  We write $q$ also for the quadralinear form obtained from $q$ and define a symmetric trilinear map $t \!: V \times V \times V \ra V$ (i.e., a linear map $S^3 V \ra V$) implicitly via
   \[
   q(x, y, z, w) = b(x, t(y, z, w)) \quad \text{for $x, y, z, w \in V$.}
   \]
   
Since $q$  has degree at least 3 and is preserved by the infinite group $G(k)$, one knows that the  hypersurface  in $\mathbb{P}(V)$ defined by the equation $q = 0$ is singular, see e.g.\ \cite[\S6]{OS}.    Consider now the singular locus $S$; it is a proper subvariety, consisting of the lines $kx$ such that $t(x,x,x)$ is zero.  We can ask if $S$ is itself smooth.  The rows of the Jacobian matrix at $kx$ are of the form $t(x, x, y)$ as $y$ varies over a basis of $V$.  We say that $x$ is \emph{is rank one} if it is nonzero and the Jacobian matrix has rank $\le 1$, i.e., the image of the linear map $y \mapsto t(x,x,y)$ has dimension $\le 1$.

\begin{eg} \label{E7.1}
The highest weight vector $\hw \in V$ is rank one.  Indeed, let $y \in V$ be a weight vector and write its weight as $\omega - \alpha$ where $\alpha$ is a sum of positive roots.  If $t(\hw, \hw, y)$ is not zero it is a weight vector with weight $3\omega - \alpha$, necessarily equal to $\omega - \alpha'$ for some sum of positive roots $\alpha'$.
However, the lowest weight of $V$ is $w_0 \omega$ where $w_0$ is the longest element of the Weyl group of $G$, so $\alpha$ is at most $\omega - w_0 \omega$.  Putting these observations together with the fact that $w_0 \omega = -\omega$, we have:
\[
2\omega = \alpha - \alpha' \le \alpha \le 2\omega.
\]
Hence $\alpha = 2\omega$ and $\alpha' = 0$.  In particular, $t(\hw, \hw, V)$ is contained in the span of $\hw$.
\end{eg}

Note the contrast with the situation for $E_6$.  In that case, the group acts transitively on the singular locus, hence the singular locus is smooth.

\begin{rmk}
Bringing to bear the explicit formula for $t(x, x, y)$ from \cite[p.~88]{Brown:E7} and representatives of the orbits for the group action from \cite[Th.~29]{Krut:E7}, one can calculate that the Jacobian of $S$ at a point $kx$ has rank 1 or 12.
\end{rmk}

We claim that the $\delta$-spaces are inner ideals for all $\delta \in \D$.
An \emph{inner ideal} is a subspace $X$ of $V$ such that $t(X, V, X)$ is contained in $X$.  It suffices to check that $V_\delta$ is an inner ideal.  The lowest weight for the action of $G$ on $V_\delta$ is $\omega_7 - \alpha$, where $\alpha$ is a sum of simple roots in the $\delta$-component, and the lowest weight for $V$ is $-\omega_7$.  Therefore, every weight of $t(V_\delta, V_\delta, V)$ is at least $\omega_7 - 2 \alpha$.  But every weight of $V$ that differs from $\omega_7$ by a sum of simple roots in the $\delta$-component belongs to $V_\delta$ by Lemma \ref{not.comp}.  We have proved that $V_\delta$ is an inner ideal.

By Prop.~\ref{incident.1}.1, the $\alpha_i$-spaces consist of rank one elements except possibly for the $\alpha_1$-spaces.  We say that an inner ideal is rank one if it consists of rank one elements.

\subsection*{Transitivity}
The group $G$ acts transitively on the 1-dimensional subspaces of $V$ spanned by rank one elements by \cite[6.2, 7.7]{Ferr:strict}, that is, the $\alpha_7$-spaces are the 1-dimensional rank one inner ideals.

The roots of $E_7$ not involving $\alpha_7$ form a root system of type $E_6$; we write $E_6$ for the corresponding subgroup of $G$.  Restricting the weights of $V$ to $E_6$, we find that $V$ is---as an $E_6$-module---a direct sum of the standard representation of $E_6$, its dual, and two 1-dimensional trivial representations.  By \cite[6.12]{G:struct}, every $d$-dimensional rank one inner ideal is in the $G$-orbit of one that is a direct sum of one of the 1-dimensional representation and a singular inner ideal of dimension $d - 1$ in the standard representation of $E_6$.  By the transitivity results from \S\ref{E6}, for $i = 2, 4, 5, 6$ the $\alpha_i$-spaces are the rank one inner ideals of the dimensions specified in \eqref{E7.dynk}.  Similarly, the $\alpha_3$-spaces are the maximal rank one inner ideals of dimension 6.

The group $G$ acts transitively on the collection of 12-dimensional inner ideals by \cite[6.15]{G:struct}, hence the $\alpha_1$-spaces are the 12-dimensional inner ideals.

\medskip

For the sake of brevity, we omit the ``incidence'' portion of this example.

\section{Loose ends} 

To complete our discussion of the concrete realization of the geometries, we address some loose ends.  Above, we have skipped the groups of type $B$ and $C$;  the reader should have no trouble filling them in from the previous examples.

\begin{borel}{Example: type $F_4$}
We now sketch the case where $G$ is of type $F_4$.  We find the following diagram:
\[
\begin{picture}(6,1)
    \multiput(1,0.5)(1,0){4}{\circle*{\darkradG}}

    \put(1,0.5){\line(1,0){1}}
    \put(3,0.5){\line(1,0){1}}
    \put(2,0.55){\line(1,0){1}}
    \put(2,0.45){\line(1,0){1}}

    \put(2.4, .35){\makebox(0.2,0.3)[s]{$<$}}

    \put(1,0.6){\makebox(0,0.4)[b]{$1$}}
    \put(2,0.6){\makebox(0,0.4)[b]{$2$}}
    \put(3,0.6){\makebox(0,0.4)[b]{$3$}}
    \put(4,0.6){\makebox(0,0.4)[b]{$6$}}

    \put(1,0.1){\makebox(0,0.4)[b]{$\beta = \alpha_4$}}
    \put(2,0.1){\makebox(0,0.4)[b]{$\alpha_3$}}
    \put(3,0.1){\makebox(0,0.4)[b]{$\alpha_2$}}
    \put(4,0.1){\makebox(0,0.4)[b]{$\alpha_1$}}

\end{picture}
\]
Klimyk's provides a $G$-invariant bilinear product on $V$, which we denote by $\#$.  
(We saw the objects $G$, $V$, $\#$ in \S\ref{E6}, where they were known as $F_4$, $V_0$, $\#$.)
The usual argument shows that the product is identically zero on the $\alpha_4$-spaces, and it is zero on the $\alpha_3$- and $\alpha_2$-spaces by Prop.~\ref{incident.1}.1.

The lowest weight of $V_{\alpha_1}$ is 
   \[
   \omega_4 - (2 \alpha_4 + 2\alpha_3 + \alpha_2) = \alpha_1 + \alpha_2 + \alpha_3
   \]
by  \ref{spindle}.  For $x, y$ weight vectors in $V_{\alpha_1}$, the product $x \# y$ has weight at least $2 \alpha_1 + 2\alpha_2 + 2\alpha_3$.  In particular, that character is not a weight of $V$ since the $\alpha_1$-coordinate of $\omega_4$ is 1.  Therefore, the product $\#$ is identically zero on the $\alpha_1$-spaces also.  Freudenthal calls $\alpha_1$-spaces ``symplecta'', since they are associated with the standard representation of a group of type $C_3$.

The group $G$ acts transitively on the $d$-dimensional subspaces of $V$ on which the product $\#$ is identically zero for $d = 1$ (the $\alpha_4$-spaces) by \cite[28.27]{Frd:E7.8}  or \cite[8.6]{Asch:E6}, $d = 2, 3$ by \cite[9.5, 9.8]{Asch:E6}, and for $d = 6$ (the $\alpha_1$-spaces) by \cite[28.22]{Frd:E7.8}; alternatively, all four cases are treated by \cite[Th.~2]{Racine}.  For all of the objects in $\G_V$, incidence is the same as inclusion.  Aside from $\alpha_1$-spaces, this is Prop.~\ref{incident.1}.2.  For cases involving an $\alpha_1$-space, one can adapt the proof of \ref{incident.1}.2, using the fact that a group of type $C_3$ acts transitively on $d$-dimensional isotropic subspaces of the standard representation for $d = 1, 2, 3$.

The space $V$ may be interpreted as the trace zero elements in an Albert algebra, i.e., a 27-dimensional simple exceptional Jordan algebra; this identifies $G$ with the group of automorphisms of the algebra.  (The papers \cite{Frd:E7.8} and \cite{Racine} cited above use this viewpoint.)  An element $x$ in such an Albert algebra has square zero if and only if the trace of $x$ is zero (i.e., $x$ is in $V$) and $x\#x = 0$, as can be easily seen using the ``sharp expression'' and ``spur formula'' from \cite{McC}.  Therefore, the objects of $\G_V$ are precisely the subspaces in the Albert algebra of dimensions 1, 2, 3, and 6 on which the multiplication is identically zero.
\end{borel}

\begin{borel}{Example: type $G_2$} \label{G2}
The geometry for a group $G$ of type $G_2$ is similar to that for type $F_4$, but everything is easier.  The dimensions are summarized in the following diagram:
\[
\begin{picture}(4,1)
    \put(1,.5){\line(1,0){1}}
    \put(1,.45){\line(1,0){1}}
    \put(1,.55){\line(1,0){1}}

    \multiput(1,.5)(1,0){2}{\circle*{\darkradG}}

    \put(1.4, .35){\makebox(0.2,0.3)[s]{$<$}}
    
    \put(1,.6){\makebox(0,0.4)[b]{1}}    
    \put(2,.6){\makebox(0,0.4)[b]{2}}

    \put(1,0.1){\makebox(0,0.4)[b]{$\beta = \alpha_1$}}
    \put(2,0.1){\makebox(0,0.4)[b]{$\alpha_2$}}
\end{picture}
\]
As in the $F_4$ case, there is a $G$-invariant bilinear product on $V$, which we denote by $\#$.  The usual argument shows that the product is identically zero on the $\alpha_1$-spaces, hence also on the $\alpha_2$-spaces by Prop.~\ref{incident.1}.

The vector space $V$ may be viewed as the trace zero elements in the split octonion algebra; this identifies $G$ with the group of automorphisms of the algebra \cite{Sp:ex}.  As in the $F_4$ case, an element $x$ in the split octonion algebra has square zero if and only if $x$ is in $V$ and $x \# x$ is zero.  It is easy to prove that $G$ acts transitively on the 1-dimensional subspaces of $V$ on which the multiplication is zero using \cite[1.7.3]{Sp:ex}.  The Cayley-Dickson process gives an explicit description of the octonion algebra, which one can use to prove that $G$ acts transitively on the 2-dimensional subspaces on which the multiplication is zero.  (This essentially solves Problem 23.54 in \cite{FH}, cf.~\ref{homo}.)  Aschbacher gives a different proof in \cite{Asch:G2}.

In summary, the $\alpha_i$-spaces are the $i$-dimensional subspaces of $V$ on which the product $\#$ is zero.  Incidence is the same as inclusion by Prop.~\ref{incident.1}.2.
\end{borel}

\begin{borel}{Example: type $E_8$}\label{E8}
The recipe in \S\ref{realize.2} for giving a concrete realization of the geometry associated with a  group has been very effective with the examples considered so far.  But what of the least familiar case, where $G$ has type $E_8$?  The recipe still works, of course, and for each fundamental representation $V$, it is easy to write down the dimension of the $\delta$-spaces for each $\delta \in \D$.   This is already interesting. 
But a problem occurs when we attempt to describe the algebraic properties that characterize the $\delta$-spaces.

For example, the smallest fundamental representation $V$ of $G$ is the adjoint representation, with highest weight $\omega_8$ and dimension $248$.  Representation theory does not provide any obvious additional structure on $V$, e.g., there is no $G$-invariant quintic form.  Therefore, the only description of the $\delta$-spaces that suggests itself is in terms of the Lie algebra structure.

The next smallest fundamental representation $V$ has highest weight $\omega_1$ and dimension 3875.  This representation has $G$-invariant bilinear and cubic forms, each determined uniquely up to a nonzero scalar multiple.  Thus there is also a $G$-invariant commutative product on $V$.  Unfortunately, one cannot simply translate the analysis in \S\ref{E6} to this case.  For example, the proof of Lemma \ref{nonzero} does not translate because the weights of $V$ are not all one orbit under the Weyl group and some of the weights occur with multiplicity greater than one.


%


\end{borel}

\begin{borel}{Projective homogeneous varieties} \label{homo}
The above examples can all be viewed from the perspective of projective homogeneous varieties, i.e., projective varieties $Y$ on which $G$ acts transitively.  We maintain our assumptions that $G$ is split simply connected and $V$ is a fundamental irreducible representation as in \S\ref{realize.2}.

There is a bijection between subsets of $\D$ and isomorphism classes of projective homogeneous varieties given by sending $S \subseteq \D$ to $Y_S := G/P_S$.  For example, $Y_\emptyset$ is a point because $P_\emptyset$ is all of $G$.

For $S$ a singleton, say $\{ \delta \}$, the $k$-points of $Y_S$ are the $\delta$-spaces in $V$.  Indeed, the $\delta$-spaces are defined to be the orbit of $V_\delta$ in the appropriate Grassmannian, and $V_\delta$ has stabilizer $P_\delta$.

\begin{eg}
For $G$ of type $B$ or $D$, the variety $Y_{\alpha_1}$ is a conic.  The other $Y_\delta$'s are families of linear subspaces of the conic.
\end{eg}

For an arbitrary subset $S \subseteq \D$, a \emph{flag of type $S$} is a collection of pairwise incident subspaces $\{ X_s \mid s \in S \}$ where $X_s$ is an $s$-space.  The flags of the extreme type $\D$ are called \emph{chambers}.  We call $\{ V_s \mid s \in S \}$ the \emph{standard flag of type $S$}.  (What we call the standard chamber is traditionally called the ``fundamental chamber''.)  We need the following consequence of the fact that $\G_V$ is a building:

\begin{prop} \emph{\cite[3.16]{Ti:BN}} Every flag of type $S$ in $\G_V$ is contained in a chamber and is in the $G$-orbit of the standard flag of type $S$.
\end{prop}

In particular, $G$ acts transitively on the collection of flags of type $S$.  The stabilizer of the standard flag is the intersection  $\cap_{s \in S} P_s$, which is $P_S$ \cite[IV.2.5, Th.~3c]{Bou:g4}.  Hence the $k$-points of $Y_S$ are the flags of type $S$ in $\G_V$.

We view the examples in the preceding sections as giving explicit descriptions of the geometry $\G_V$ as well as the projective homogeneous varieties under split groups $G$.
When $G$ is not split, the situation is somewhat more complicated.  The absolute Galois group of $k$ acts on the Dynkin diagram $\D$, and there is a bijection between Galois-invariant subsets $S$ of $\D$ and projective homogeneous varieties defined over $k$.  The description in the split case can be altered to give a description in the general case.  For groups of ``inner type'' $A_n$, one finds the generalized Severi-Brauer varieties as in \cite[\S1]{KMRT}.  (Note that these varieties may have no $k$-points.)
For groups of type $D_n$, see \cite{G:flag} ($n = 4$) and \cite[p.~183]{MPW2} ($n \ne 4$).
\end{borel}

\section{Outer automorphisms} \label{outer}

Let $\G_V$ be a geometry defined from an irreducible representation $V$ of $G$ by the recipe in \S\ref{realize}.
Every automorphism $\phi$ of $G$ permutes the parabolic subgroups, hence induces an automorphism of Tits's geometry $\G_P$.  Further, $\phi$ induces an automorphism of the concrete geometry $\G_V$ via the isomorphism between  $\G_P$ and $\G_V$ from \S\ref{realize}.

Every $g \in G(k)$ defines an automorphism of $G$ by sending $h \mapsto g h g^{-1}$.  Such an automorphism is called \emph{inner}.  It preserves types of objects and sends $X \in \G_V$ to $gX$.   In classical projective geometry, such an automorphism is called a \emph{collineation}.

On the other hand, some groups have automorphisms that are not of this type; such automorphisms are called \emph{outer}.  They have a more interesting action on the geometry $\G_V$ in that they do not preserve the types of objects.  In classical projective geometry, they are called \emph{correlations}. 
As an example, the map $g \mapsto (g^{-1})^t$ is an automorphism of $SL_3$, and it is outer because it does not fix the center elementwise.  We will see in Example \ref{An.dual} below that the induced map $\psi$ is the polarity with respect to a certain conic.  

Generally speaking, the existence of an outer automorphism of $G$ implies a principle of duality (for $D_4$, triality) in the geometry $\G_V$.  For $SL_3$---equivalently, $\mathbb{P}^2$---it takes the following form \cite[2.3]{Coxeter}: ``every definition remains significant, and every theorem remains true, when we interchange \emph{point} and \emph{line}, \emph{join} and \emph{intersection}.''  See \cite[p.~155]{Weiss} for an analogous statement of the principle of triality.

Let $\phi$ be an automorphism of $G$, and let $\sub(V)$ denote the collection of subspaces of $V$.
We want an efficient way to check if a given function $\psi \!: \G_V \ra \sub(V)$ is the automorphism of the geometry $\G_V$ induced by $\phi$.

\begin{thm} \label{carr}
If 
\begin{enumerate}
\item $\psi(g X) = \phi(g) \psi(X)$ for every $X \in \G_V$ and $g \in G$ and
\item there is a chamber $\{ V_i \mid 1 \le i \le n \}$ such that $\{ \psi(V_i) \mid 1 \le i \le n \}$ is also a chamber,
\end{enumerate}
then $\psi$ is the automorphism of the geometry $\G_V$ induced by the automorphism $\phi$ of $G$.
\end{thm}

[The term ``chamber'' was defined in \ref{homo}.]

In the examples below, we will specify a function $\psi$ and prove that it satisfies the hypotheses of the theorem.  We will use all of the freedom implicit in our hypotheses on $\psi$.  In the $E_6$ example, it will be clear that $\psi(X)$ is a subspace of $V$ for $X \in \G_V$, but it will not be obvious that $\psi(X)$ is in $\G_V$.  In the triality example, we will only define $\psi$ on $\G_V$ and not for an arbitrary subspace of $V$.  And in the type $A$ example, $\psi$ will not preserve the standard chamber.

\begin{proof}[Proof of Theorem \ref{carr}]
Let $X$ be an object in $\G_V$.
We first claim that $\psi(X)$ is also in $\G_V$.  We find a chamber $\{ X_i \mid 1 \le i \le n \}$ containing $X$ such that $X_i$ is of type $\alpha_i$.  This chamber is conjugate to the chamber $\{ V_i \mid 1 \le i \le n \}$ from (2), i.e., there is some $g \in G$ such that $gV_i = X_i$ for every $i$.  Therefore $\psi(X_i) = \phi(g) \psi(V_i)$, and $\psi(X_i)$ is an object in the geometry for all $i$.

Let $P$ be the stabilizer of $X$ in $G$.  For $g \in \phi(P)$, we have
   \[
   g \psi(X) = \psi(\phi^{-1}(g) X) = \psi(X) \quad \text{by (1),}
   \]
hence $\phi(P)$ is contained in the stabilizer of $\psi(X)$.  But $\psi(X)$ is an object in $\G_V$, hence by definition it is a nonzero, proper subspace of $V$.  In particular, its stabilizer is a proper subgroup of $G$.  Since $P$ is a maximal proper subgroup of $G$, so is $\phi(P)$, hence $\phi(P)$ is the stabilizer of $\psi(X)$.  This proves that the diagram
   \[
   \begin{CD}
   \G_P @>\phi>> \G_P \\
   @AAA @AAA \\
   \G_V @>\psi>> \G_V
   \end{CD}
   \]
commutes, where the vertical arrows send a subspace of $V$ to its stabilizer in $G$.  Since the vertical arrows are bijections (see \S\ref{realize}), $\psi$ is also a bijection.  
Moreover, $\psi$ respects the notion of incidence in $\G_V$, because that relation is the one transported from $\G_P$ by the vertical isomorphisms. We have proved that $\psi$ is an automorphism of $\G_V$, and the commutativity of the diagram shows that it is the one induced by $\phi$.
\end{proof}

\begin{eg}[type $A$: projective duality] \label{An.dual} 
Let $G$ be $SL_n$ acting on $k^n$, and let $\phi$ be the automorphism $g \mapsto (g^{-1})^t$.  For a subspace $X$ of $k^n$, we define $\psi(X)$ to be the orthogonal complement of $X$ with respect to the dot product defined by $x \cdot y := x^t \, y$.

Viewed algebraically, the dot product identifies $k^n$ with the dual vector space $(k^n)^*$.  This identification pairs $\psi(X)$ with the collection of linear forms vanishing on $X$.

Viewed geometrically, the map $\psi$ is precisely the correspondence between points and hyperplanes giving projective duality in $\mathbb{P}^{n-1}$ described in \cite{Pedoe} and \cite[11.8]{Coxeter}.  It  interchanges a point $[ a_1 : a_2 : \ldots : a_n ]$ in homogeneous coordinates with the hyperplane consisting of solutions to the equation $\sum a_i x_i = 0$.  For $n = 3$, it is the polarity with fundamental conic $x_1^2 + x_2^2 + x_3^2 = 0$, cf.~\cite[\S98]{VY}.

We now check that $\psi$ satisfies the conditions of Theorem \ref{carr}.  
The dot product is compatible with the automorphism $\phi$ in the sense that
   \[
   x \cdot y = (gx) \cdot (\phi(g)y)   \quad \text{for $g \in SL_n$ and $x, y \in k^n$.}
   \]
Thus, a vector $y$ is in $\psi(gX)$ if and only if $X \cdot (\phi(g)^{-1} y) = 0$, i.e., if and only if $y$ is in $\phi(g) \psi(X)$.  Thus $\psi$ satisfies (1).

Consider the collection $\{ V_1, V_2, \ldots, V_{n - 1} \}$ of subspaces such that $V_i$ consists of the vectors whose bottom $n - i$ coordinates are zero.  Each $V_i$ is stabilized by the upper triangular matrices---which make up a Borel subgroup---so this collection is a chamber.  
Applying $\psi$, we find that $\psi(V_i)$ is the space of vectors whose top $n - i$ coordinates are zero.
Each $\psi(V_i)$ is stabilized by the lower triangular matrices.  That is, $\{ \psi(V_1), \psi(V_2), \ldots, \psi(V_{n-1}) \}$ is also a chamber and $\psi$ satisfies (2). Hence $\psi$ is the automorphism of $\G_V$ corresponding to the automorphism $\phi$ of $SL_n$.
\end{eg}

\section{Example: type $D$ (orthogonal duality)} \label{Dn.dual}
Let $G$ be a group as in \S\ref{Dn}, constructed from the root system of type $D_n$ for some $n \ge 4$; it is traditionally denoted by $\Spin_{2n}$.  Its Dynkin diagram has an automorphism of order 2 given by interchanging the roots $\alpha_{n-1}$ and $\alpha_{n}$.  Here we will construct the corresponding automorphism $\psi$ of the geometry $\G_V$.

\begin{borel*} \label{Dn.b}
We can draw a Hasse diagram for the weights of $V$ as we did for $E_6$ in Figure \ref{weights}.   The case $n = 4$ is shown in Figure \ref{D4.weights} below. 
Figure \ref{Dn.weights} shows the diagram in the general case, rotated counterclockwise by $90^\circ$ for space considerations. 
\renewcommand{\thefigure}{\theequation}%
\refstepcounter{equation}\label{Dn.weights}
\begin{figure}[tbh] 
\begin{picture}(10.05,1.5)
    \put(.3,.75){\line(1,0){2.1}}
    \put(3.3,.75){\line(1,0){1.2}}
    \put(5.775,.75){\line(1,0){1.2}}
    \put(7.875,.75){\line(1,0){2.1}}
    
    \put(4.5,.75){\line(1,1){.621}}
    \put(4.5,.75){\line(1,-1){.621}}
    
    \put(5.1375,.1125){\line(1,1){.621}}
    \put(5.1375,1.3875){\line(1,-1){.621}}
    
    \multiput(.3,.75)(.9,0){3}{\circle*{\darkradG}}
    \multiput(3.6,.75)(.9,0){2}{\circle*{\darkradG}}
    
    \multiput(8.175,.75)(.9,0){3}{\circle*{\darkradG}}
    \multiput(5.775,.75)(.9,0){2}{\circle*{\darkradG}}
    
    \put(5.1375,1.3875){\circle*{\darkradG}}
    \put(5.1375,.1125){\circle*{\darkradG}}
    
    \multiput(2.625,.75)(.225,0){3}{\circle*{0.01}}
    \multiput(7.2,.75)(.225,0){3}{\circle*{0.01}}
    
    \put(0.75,0.35){\makebox(0,0.4)[b]{1}}    
    \put(1.65,0.35){\makebox(0,0.4)[b]{2}}
    
    \put(9.525,0.35){\makebox(0,0.4)[b]{1}}    
    \put(8.625,0.35){\makebox(0,0.4)[b]{2}}
    
    \put(4.05,0.35){\makebox(0,0.4)[b]{$n - 2$}}    
    \put(6.225,0.35){\makebox(0,0.4)[b]{$n - 2$}}
    
    \put(5.85, 1.125){\makebox(0,0.4)[b]{$n - 1$}}
    \put(4.7, 1.125){\makebox(0,0.4)[b]{$n$}}
    \put(5.4, .075){\makebox(0,0.4)[b]{$n$}}
    \put(4.4, .075){\makebox(0,0.4)[b]{$n - 1$}}
\end{picture}
\caption{Hasse diagram of weights of the standard representation of $D_n$ from \cite[Fig.~4]{PSV}.  Larger weights are on the left.}
\end{figure}
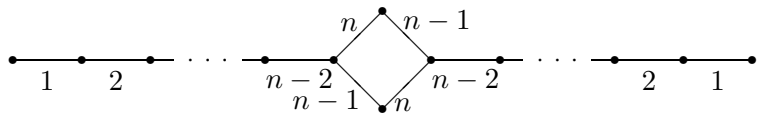

Fix nonzero vectors $e_1, e_2, \ldots, e_n$ in $V$ such that $e_i$ has weight
\[
\e_i := \omega_1 - \sum_{j = 1}^{i - 1} \alpha_j.
\]
The longest element of the Weyl group is $-1$; it is the unique automorphism of the diagram of order 2 that stabilizes none of the weights.  The weights $\e_1$ through $\e_{n-1}$ are those in the string on the left of the diagram and $\e_n$ is the bottom weight in the middle square.  The other weights of $V$ are of the form $-\e_i$ for some $i$.  Let $f_i$ be a nonzero vector of weight $-\e_i$, so the vectors $e_1, e_2, \ldots, e_n, f_1, f_2, \ldots, f_n$ are a basis of $V$.

Let $b$ be the $G$-invariant symmetric bilinear form on $V$ as in \S\ref{Dn}.  Clearly, since $\e_i$ is not $-\e_j$ for any pair $i, j$, the subspace of $V$ spanned by the $e_i$'s (respectively, by the $f_i$'s) is isotropic, i.e., $b(e_i, e_j) = b(f_i, f_j) = 0$ for all $i, j$.   Also, $b(e_i, f_j)$ is nonzero if and only if $i = j$.  By scaling the $f_j$'s, we may assume that $b(e_i, f_j) = \delta_{ij}$ (Kronecker delta).  (We have now obtained the description of $G$ and $SO(b)$ given in \cite[\S{V.7}]{Brown}.)  The construction in \S\ref{Dn} gives:
   \begin{gather}
   V_{\alpha_i} := \kspan \{ e_1, e_2, \ldots, e_i \}\ \text{for $i \le n - 2$}, \notag \\
   V_{\alpha_{n-1}} = \kspan \{ e_1, e_2, \ldots, e_n \}, \text{and}   \label{Dn.chamber} \\
    V_{\alpha_n} = \kspan \{ e_1, e_2, \ldots, e_{n-1}, f_n \}. \notag
   \end{gather}
\end{borel*}

Let $s$ denote the matrix in $GL(V)$ that fixes $e_i$ and $f_i$ for $1 \le i < n$ and interchanges $e_n$ and $f_n$.  It leaves $b$ invariant, and the map $\phi \!: SO(b) \ra SO(b)$ defined by $\phi(g) = sgs^{-1}$ is an automorphism of order 2.  There is a unique lift of $\phi$ to an automorphism of $G$ \cite[2.24(i)]{BoTi:C}, which we also denote by $\phi$.  The description of the root subgroups in $SO(b)$ in \cite[23.4]{Borel} shows that $\phi$ is, in fact, the automorphism of $G$ induced by the automorphism of the Dynkin diagram that interchanges $\alpha_{n-1}$ and $\alpha_n$.

For each subspace $X$ of $V$, put $\psi(X) := sX$.  It is a triviality that $\psi$ satisfies condition (1) of \ref{carr} and that the standard chamber exhibited in \eqref{Dn.chamber} is permuted by $\psi$, hence that $\psi$ satisfies condition (2).
That is, $\psi$ is the automorphism of $\G_V$ induced by the automorphism $\phi$ of $G$ and $SO(b)$.

\section{Example: type $D_4$ (triality)} \label{D4}
Continue the notation of the preceding section, \S\ref{Dn.dual}, except suppose now that $n = 4$.  The Dynkin diagram of $G$ looks like
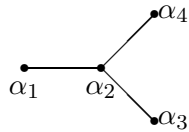
\[
\begin{picture}(4, 2)
    \put(1,1){\line(1,0){1}}

    \put(2,1){\line(1,-1){0.7}}
    \put(2,1){\line(1,1){0.7}}

    \multiput(1,1)(1,0){2}{\circle*{\darkradG}}
    
    \put(1,0.6){\makebox(0,0.4)[b]{$\alpha_1$}}  
    \put(2,0.6){\makebox(0,0.4)[b]{$\alpha_2$}}  

    \put(2.7,0.3){\circle*{\darkradG}}
    \put(2.7,1.7){\circle*{\darkradG}}

    \put(2.95,0.2){\makebox(0,0.4)[b]{$\alpha_3$}}  
    \put(2.95,1.6){\makebox(0,0.4)[b]{$\alpha_4$}}

\end{picture}
\]
Let $\phi$ be the automorphism of order 3 that permutes the arms counterclockwise.
We will now describe explicitly the corresponding automorphism $\psi$ of the geometry $\G_V$.  

Let $\rho_0$ be the representation of $G$ on $V$ with highest weight $\omega_1$.  For $i = 1, 2$, we set $\rho_i := \rho_0 \phi^{-i}$; it is a representation of $G$ on $V$.  The highest weight of $\rho_1$ is $\phi(\omega_1) = \omega_3$, and the highest weight of $\rho_2$ is $\phi^2(\omega_1) = \omega_4$.

The weights of $V$ with respect to $\rho_0$ are listed in Figure \ref{D4.weights}; such a diagram is easily constructed as in \ref{min}.  The weights relative to $\rho_1$ and $\rho_2$ are the same except with $\phi$ or $\phi^2$ applied, respectively.   Fix a basis $e_i$, $f_j$ for $V$ and a symmetric bilinear form $b$ as in \ref{Dn.b}.
The image $\rho_i(G)$ of $G$ in $GL(V)$ is the same for all $i$, so $b$ is $\rho_i(G)$-invariant for all $i$.
\refstepcounter{equation} \label{D4.weights}
\begin{figure}[th]
\begin{minipage}{\textwidth}
\[
\xymatrix@R=8pt{
& (1, 0, 0, 0) = \omega_1 \ar@{-}[d]^1 & \\
& (-1, 1, 0, 0) \ar@{-}[d]^2 & \\
& (0, -1, 1, 1) \ar@{-}[dl]_3 \ar@{-}[dr]^4 & \\
(0, 0, -1, 1) \ar@{-}[dr]_4 & & (0, 0, 1, -1) \ar@{-}[dl]^3 \\
& (0, 1, -1, -1) \ar@{-}[d]^2 & \\
& (1, -1, 0, 0) \ar@{-}[d]^1 & \\
& (-1, 0, 0, 0) = -\omega_1
}
\]
\end{minipage}
\caption{Hasse diagram of the weights of $V$ relative to $\rho_0$.}
\end{figure}
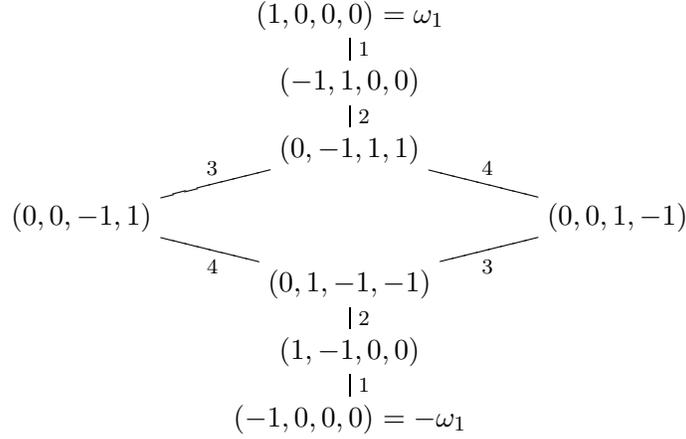

We claim that there is a nonzero linear map $t \!: V \ot V \ot V \ra k$ that is $G$-invariant in the sense that
\begin{equation} \label{trial.1}
t(\rho_0(g) x_0, \rho_1(g) x_1, \rho_2(g) x_2) = t(x_0, x_1, x_2)
\end{equation}
and is unique up to multiplication by a nonzero scalar.  To see this, note that a linear map satisfying \eqref{trial.1} is a $G$-fixed element of $\rho_0^* \ot \rho_1^* \ot \rho_2^*$, where $\rho_i^*$ denotes the representation dual to $\rho_i$---which in this case is just $\rho_i$.  We compute the formal character of this tensor product by multiplying the formal characters of the $\rho_i$ \cite[22.5B]{Hum:LA}, recover the decomposition of $\rho_0 \ot \rho_1 \ot \rho_2$ as a direct sum of irreducibles, and look for 1-dimensional summands as in \ref{repth}.  (The software package LiE uses the more efficient Brauer-Klimyk formula described in \cite[Exercise 24.9]{Hum:LA}.)  We find a unique 1-dimensional summand; existence and uniqueness of $t$ follows.

We can prove that $t$ is nonzero for some specific arguments.

\begin{lem} \label{trial.nonzero}
Let $x_i \in V$ be a nonzero vector of weight $\mu_i$ relative to $\rho_i$ for $i = 0, 1, 2$.  We have: $t(x_0, x_1, x_2)$ is nonzero if and only if $\mu_0 + \mu_1 + \mu_2 = 0$.
\end{lem}

\begin{proof}
``Only if'' is clear, so we prove ``if''.  Suppose that $\sum \mu_i =  0$.  Since $\rho_2(G)$ acts transitively on the weights of $V$ relative to $\rho_2$, we may assume that $\mu_2$ is $\omega_4$.  By the argument in the proof of Lemma \ref{orbits}, the subgroup of the Weyl group fixing $\omega_4$ has two orbits on the weights of $V$ relative to $\rho_1$, with representatives $\pm \omega_3$.  Since $-\omega_4 - \omega_3$ is not a weight of $V$ relative to $\rho_0$, we must have $\mu_1 = -\omega_4 + \omega_3$.  We have just proved that the Weyl group acts transitively on the triples $\mu_0, \mu_1, \mu_2$ such that $\sum \mu_i = 0$.  As in the proof of Lemma \ref{nonzero}, it follows that $t(x_0, x_1, x_2)$ is nonzero.
\end{proof}

Moreover, $t$ is invariant under cyclic permutations.

\begin{lem} \label{cyclic}
The value of $t$ is unchanged if its arguments are permuted cyclically.
\end{lem}

\begin{proof}
Consider the linear map $d \!: V \ot V \ot V \ra k$ defined by 
  \[
  d(x_0, x_1, x_2) := t(x_0, x_1, x_2) - t(x_1, x_2, x_0);
  \]
it is $G$-equivariant because $\rho_{i+1} = \rho_i \phi^{-1}$ for all $i$.
By the uniqueness of $t$, the map $d$ must be a scalar multiple of $t$.  The vector $e_4 \in V$ is nonzero of weight $-\omega_3 + \omega_4$ relative to $\rho_0$.  Then $t(e_4, e_4, e_4)$ is not zero by the previous lemma, yet $d(e_4, e_4, e_4)$ is zero.  Therefore, $d$ is identically zero.
\end{proof}

We now define products $*_i$ on $V$ for $i = 0, 1, 2$ implicitly via
\[
t(x_0, x_1, x_2) = b(x_i, x_{i+1} *_i x_{i+2}).
\]
By Lemma \ref{cyclic}, all three products agree, so we write simply $*$.  Because $t$ and $b$ are $G$-equivariant, in the sense that
\begin{equation} \label{trial.2}
(\rho_i(g) x) * (\rho_{i+1}(g) y) = \rho_{i+2}(g) (x * y).
\end{equation}
This allows us to compute the multiplication, at least up to a scalar factor.  Let $x_i \in V$ be nonzero with weight $\mu_i$ relative to $\rho_i$ for $i = 1, 2$.  It follows from Lemma \ref{trial.nonzero} that $x_1 * x_2$ is nonzero if and only if $\mu_1 + \mu_2$ is a weight of $V$ relative to $\rho_0$, in which case $x_1 * x_2$ has weight $\mu_1 + \mu_2$.  We summarize these computations in the table below, where the entry in the row $x_1$ and column $x_2$ is ``$\cdot$'' if $x_1* x_2$ is zero and, for example, $e_3$ if $x_1 * x_2$ is a nonzero scalar multiple of $e_3$.  The left column lists the weight of $x_1$ for the reader's convenience; we omit the weight of $x_2$ due to space considerations.  Since the product is $G$-equivariant and the weights of $\rho_i$ are preserved under multiplication by $-1$, one needs only to compute the first four columns of entries; the remaining four columns can be filled in by symmetry.
\begin{equation}\label{multtable}
\begin{array}{ccc|cccc|cccc} 
&   &     &\multicolumn{8}{c}{x_2} \\
&  &     &  e_1&  e_2&  e_3&  e_4&  f_4&  f_3&  f_2&  f_1\\ \hline
&(0,0,1,0)    &e_1&\cdot&\cdot&\cdot& e_1&\cdot& e_2&  e_3&f_4\\
&(0, 1,-1,0)   &e_2&\cdot&\cdot&  e_1&\cdot& e_2&\cdot& e_4& f_3\\
 &(1,-1,0,1)  &e_3&\cdot& e_1&\cdot&\cdot& e_3& e_4&\cdot&  f_2\\
x_1&(1,0,0,-1)&e_4&e_1& \cdot& \cdot& \cdot&f_4&f_3&f_2& \cdot \\ \hline
 &(-1,0,0,1)  &f_4& \cdot&e_2&e_3&e_4& \cdot& \cdot& \cdot&f_1\\
  &(-1,1,0,-1) &f_3&  e_2&\cdot& f_4& f_3&\cdot&\cdot& f_1&\cdot\\
  &(0, -1, 1, 0) &f_2& e_3& f_4&\cdot& f_2&\cdot&  f_1&\cdot&\cdot\\
  &(0,0,-1,0) &f_1& e_4&  f_3& f_2&\cdot& f_1&\cdot&\cdot&\cdot
\end{array}
\end{equation}

Although the table above is not fine enough to allow us to actually multiply two vectors in $V$, it is sufficient to describe how the objects in the geometry $\G_V$ interact with the multiplication.  Specifically, we can recover some of the results of \cite{vdBSp} without discussing octonion algebras.

\begin{prop} {\rm (Cf.~\cite[\S2]{vdBSp})} \label{vdBSp.prop}
\begin{enumerate}
\item If $X$ is an $\alpha_1$-space (a ``point''), then $V * X$ is an $\alpha_3$-space and $X * V$ is an $\alpha_4$-space.
\item If $X$ is an $\alpha_2$-space (a ``line''), then $(X * V) * X$ and $X * (V * X)$ are $\alpha_2$-spaces.
\item If $X$ is an $\alpha_3$-space (resp., an $\alpha_4$-space), then there is a unique $\alpha_1$-space $U$ such that $X = V * U$ (resp., $X = U * V$).
\end{enumerate}
\end{prop}

\begin{proof}
By \eqref{trial.2}, it suffices to check (1) for the case where $X$ is the $k$-span of $e_1$, i.e., $V_{\alpha_1}$.  In that case, (1) is clear from the multiplication table.  A similar argument handles (2).

We now prove (3) for $\alpha_3$-spaces.  As in the previous paragraph, it suffices to check the case where $X$ is the $k$-span of $e_1, e_2, e_3, e_4$, i.e., $V_{\alpha_3}$.  The multiplication table shows that $X$ is $V * e_1$, so suppose that $u \in V$ is nonzero and satisfies $V * u = X$.  Considering the product $f_1 * u \in X$, we see from the multiplication table that the coefficients of $e_2, e_3$, and $f_4$ in $u$ are all zero.  Replacing $f_1 * u$ with $f_2 * u$, etc., we conclude that $u$ is in $ke_1$.
\end{proof}

Define $\psi$ via
\begin{equation} \label{psi.t1}
\psi(ka) = a * V, \quad \psi(a * V) = V * a, \quad \text{and} \quad \psi(V * a) = ka
\end{equation}
for $a$ isotropic in $V$.  
This is well defined by the proposition. For $X$ an $\alpha_2$-space, we define
  \begin{equation} \label{psi.t2}
  \psi( X ) = X * (V * X).
  \end{equation}

Applying \eqref{trial.2}, it is easy to check that $\psi$ satisfies condition (1) of Th.~\ref{carr}.
On the other hand, we checked in the proof of Prop.~\ref{vdBSp.prop} that $\psi$ maps the standard chamber to the standard chamber, so $\psi$ also satisfies condition (2).  Thus $\psi$ is the automorphism of $\G_V$ corresponding to $\phi$.

\begin{rmks*}
Equation \eqref{psi.t1} defines $\psi^2$ on the $\alpha_1$-, $\alpha_3$-, and $\alpha_4$-spaces.  Appealing to Th.~\ref{carr}, one finds that $\psi^2(X)$ is $(X * V) * X$ for $X$ an $\alpha_2$-space.

We remark that we have recovered a multiplication of the octonions---at least approximately---entirely from first principles of representation theory.
\end{rmks*}

\section{Example: type $E_6$ (duality)} \label{E6.dual}

Let $G$ be the split simply connected group of type $E_6$ with standard representation $V$ as in \S\ref{E6}.  In this section, we will give an explicit description of the automorphism $\psi$ of the geometry $\G_V$ corresponding to the automorphism $\phi$ of the group $G$.

We define the \emph{brace product} on $V$ following \cite[p.~190]{McC}:
\[
\{ x, y, z \} := b(x,y)z + b(z,y)x - (x \# z) \# y.
\]
Note that $x$ and $z$ are interchangeable.  (See Remark \ref{brace.rmk} for comments on why we choose to work with the brace product.)

For each subspace $X$ of $V$, we set
\[
\fbox{$\psi(X) := \{ y \in V \mid \{ X, y, V \} \subseteq X \}$}
\]
From \eqref{duality}, we find that
\begin{equation} \label{duality.br}
g \{ x, y, z \} = \{ gx, \phi(g) y, g z \} \quad \text{for $g \in G$ and $x, y, z \in V$.}
\end{equation}
An argument nearly identical to the one in Example \ref{An.dual} shows that $\psi$ satisfies hypothesis (1) of Th.~\ref{carr}.  The rest of this section is spent proving that $\psi(V_{\delta}) = V_{\phi(\delta)}$ for all $\delta \in \D$, i.e., $\psi$ permutes the objects in the standard chamber.  This will show that $\psi$ satisfies hypothesis (2) of the theorem.

\begin{eg} \label{crucial.eg} \label{crucial.eg1}
Let $U$ denote the set of weights $\mu$ of $V$ such that $\phi(\omega_1 + \mu)$ is \emph{not} a weight of $V$.  Let $z$ be a nonzero vector of weight $\mu \in U$; we claim that $z$ is in $\{ \hw, \lw, V \}$, where $\lw$ is a lowest weight vector of $V$.  First observe that since $z$ does not have weight $-\omega_6$, $b(\hw, z)$ is zero.  Since $\mu$ is in $U$, $\hw \# z$ is zero and we have: $\{ \hw, \lw, z \} = b(\hw, \lw) z$.  But $b(\hw, \lw)$ is not zero because $b$ is nondegenerate.  This proves the claim.

The image of the map $x \mapsto \hw \# x$ is 10-dimensional, so the subspace of $V$ spanned by weight vectors with weights in $U$ is 17-dimensional.  Therefore, 
\[
\dim \{ \hw, \lw, V \} \ge 17.
\]
\end{eg}

\begin{borel}{Connection with Jordan theory, part II} \label{conjord2}
Let $N$ and $b$ be as in \ref{conjord}.  That is, we suppose for the moment that $k$ is algebraically closed and view $N$ and $b$ as the generic norm and trace bilinear form on an Albert algebra, respectively.  We will use three facts from the theory of Albert algebras.  The first is the \emph{5-linear identity} from \cite[p.~202]{McC}:
\[
\{ x, y, \{ z, w, u \}\} = \{ \{ x, y, z \}, w, u \} - \{ z, \{ y, x, w \}, u \} + \{ z, w, \{ x, y, u \} \}.
\]

Second, we will use the classification of the inner ideals of $V$.  A subspace $X$ is an \emph{inner ideal} if $\{ X, V, X \}$ is contained in $X$.  By \cite[\S7]{McC:inn}, 
 the proper inner ideals are the singular subspaces---which all have dimension $\le 6$---and the hyperlines.  The $\alpha_2$- and $\alpha_6$-spaces are \emph{maximal} proper inner ideals.
A straightforward application of the 5-linear identity gives: \emph{If $I$ is an inner ideal in $V$, then $\psi(I)$ is also an inner ideal.}

Finally, we will use the \emph{adjoint identity}:
\begin{equation} \label{adjoint}
x^{\#\#} = N(x)x \quad \text{for $x \in V$.}
\end{equation}
(Of course, this identity can also be proved using the representation-theoretic fact that there is a unique $G$-invariant line in $V \ot S^4(V^*)$.)

We claim that the three facts hold also in the case where $k$ is not algebraically closed.  For the adjoint identity, this is clear.  Similarly, a subspace $X$ of $V$ is an inner ideal if and only if $X \ot \kbar$ is an inner ideal in $V \ot \kbar$.  The only potentially tricky check is the claim that every nonsingular, nonzero, and proper inner ideal $X$ is a hyperline.  To see this, note that $X \ot \kbar$ is a hyperline by the case where $k$ is algebraically closed, so $X \ot \kbar$ equals $(X \ot \kbar)^\# \# (V \ot \kbar)$ by Remark \ref{hyper.quad}.  But we may rewrite this as $(X^\# \# V) \ot \kbar$.  Since $(X \ot \kbar)^\#$ is 1-dimensional and singular, so is $X^\#$.  This shows that $X$ is the hyperline $X^\# \# V$.
\end{borel}

\begin{borel}{Computation of $\psi(V_{\alpha_1})$} \label{Va1}
Linearizing the adjoint identity as in \cite[p.~496]{McC:FST}, we find the identity (McCrimmon's Equation (12)):
   \begin{equation} \label{eq12}
  (x \# y) \# (x \# z) = b(x^\#, y) z + b(x^\#, z) y + b(y \# z, x) x - x^\# \# (y \# z)
  \end{equation}
For the highest weight vector $\hw$, we have $(\hw)^\# = 0$ and \eqref{eq12} becomes
  \[
  (\hw \# y) \# (\hw \# z) = b(y \# z, \hw) \hw.
  \]
We find:
  \[
  \{ \hw, {\hw} \# z, y \} = b(\hw, {\hw} \# z) y + b(y, {\hw} \# z) \hw - b(y\#z, \hw) \hw.
  \]
Since $N(-, -, -)$ is symmetric, equation \eqref{sharp.def} shows that the first summand is zero and the second and third summands cancel.  Therefore, $\{ \hw, {\hw} \# z, y \}$ is zero and 
$\psi(k\hw)$ contains ${\hw} \# V$.

Since $\psi(k\hw)$ is an inner ideal containing the hyperline ${\hw} \# V$ and it is proper  (Example \ref{crucial.eg}), the ideal must be the hyperline.
Since $G$ acts transitively on the singular 1-dimensional subspaces and $\psi$ satisfies \ref{carr}.1, we obtain the following lemma:

\begin{lem} \label{dim1}
If $X$ is a 1-dimensional singular subspace (= an $\alpha_1$-space), then $\psi(X)$ is the hyperline $X \# V$.
\end{lem}
\end{borel}

\begin{borel}{Computation of $\psi(V_{\alpha_2})$} \label{Va2}
We now show that $\psi(V_{\alpha_2})$ is $V_{\alpha_2}$.  For $x, y$ weight vectors in $V_{\alpha_2}$ and $w$ a weight vector in $V$, we find that $\{ x, y, w \}$ has weight at least
   \[
   \la_2 + \phi(\la_2) - \omega_6 = \omega_1 - 2(\alpha_1 + \alpha_3 + \alpha_4 + \alpha_5 + \alpha_6).
   \]
But every weight of $V$ of the form $\omega_1 - (c_1 \alpha_1 + c_3 \alpha_3 + c_4 \alpha_4 + c_5 \alpha_5 + c_6 \alpha_6)$ with each $c_i$ a nonnegative integer belongs to $V_{\alpha_2}$ by Lemma \ref{not.comp}.  Hence $V_{\alpha_2}$ is contained in $\psi(V_{\alpha_2})$.

Since $\lw$ is not in $\psi(V_{\alpha_2})$ by Example \ref{crucial.eg} and $\psi(V_{\alpha_2})$ is an inner ideal, the classification of inner ideals gives that  $\psi(V_{\alpha_2})$ is precisely $V_{\alpha_2}$.
\end{borel}


\begin{eg} \label{crucial.eg2}
Let $y$ be a nonzero vector of weight $\la_2$, i.e., a lowest weight vector for $L_{\alpha_2}$ acting on $V_{\alpha_2}$.  We claim that $\{ \hw, y, V \}$ is $V_{\alpha_2}$.  Since $y$ is in $\psi(V_{\alpha_2})$ by \ref{Va2}, we need only show that $V_{\alpha_2}$ is contained in $\{ \hw, y, V \}$.  It suffices to prove that every weight of $V_{\alpha_2}$ is a weight of $\{ \hw, y, V \}$ because every weight of $V$ has multiplicity one.

Consulting Figure \ref{weights}, we see that the weights of $V_{\alpha_2}$ are symmetric in the following sense: If $\omega_1 - \alpha$ is a weight of $V_{\alpha_2}$, then $\la_2 + \phi(\alpha)$ is also a weight of $V_{\alpha_2}$.  Consequently, for every weight $\la$ of $V_{\alpha_2}$, 
  \[
  f(\la) := -\phi (\la_2 + \phi(\omega_1 - \la)) = -\phi(\la_2) + \la - \omega_1
  \]
is a weight of $V$.
For each weight $\la$ of $V_{\alpha_2}$, fix a nonzero vector $z_\la$ of weight $f(\la)$.  The vector $\{ \hw, y, z_\la \}$ has weight $\la$, and it suffices to prove that it is not zero for each $\la$.  Note that $b(\hw, y)$ is zero because $\omega_1 + \phi(\la_2)$ is not zero.

For $\la = \omega_1$, we note that $\omega_1 + f(\la) = \omega_1 - \phi(\la_2)$ has 2 as one of its entries, hence $\phi(\omega_1 + f(\la))$ is not a weight of $V$.  Therefore, ${\hw} \# z_\la$ is zero.  We have
   \[ 
   \{ \hw, y, z_\la \} = b(z_\la, y) \hw,
   \]
which is not zero because $f(\la) = -\phi(\la_2)$.

For the other five weights $\la$ of $V_{\alpha_2}$, note that $f(\la) \ne -\phi(\la_2)$, so $b(z_\la, y)$ is zero.  We claim that ${\hw} \# z_\la$ is not zero.  It has weight
$\mu := \phi(\omega_1 + f(\la)) = \phi(\la) - \la_2$.  For $\la = \la_3$,  Equation \eqref{lambdas} gives that $\mu = \la_5$, a weight of $V$.  For $\la = \la_4 = \la_3 - \alpha_3$, we find that $\mu$ is $\la_5 - \alpha_5$.  Similarly, we find that for each of the the three remaining $\lambda$'s, the weight $\mu$ is a weight of $V$.
That is, ${\hw} \# z_\la$ is nonzero.  The function $f$ was defined so that $({\hw} \# z_\la) \# y$ would have weight $\la$, hence that product is also not zero, i.e.,
  \[
  \{ \hw, y, z_\la \} = ({\hw} \# z_\la) \# y \ne 0.
  \]
We have proved that $\{ \hw, y, V \}$ is $V_{\alpha_2}$.
\end{eg}

We can now give a reasonably good description of the space $\{ x, y, V \}$ when $x$ and $y$ are nonzero weight vectors.  Note that $x$ and $y$ are necessarily singular because no weight of $V$ has 2 as a coordinate.  We say that $y$ and $x \#V$ are \emph{connected} if there is a singular vector $z \in x \# V$ such that $y$ and $z$ are ``collinear'', i.e., such that $y$ and $z$ span an $\alpha_3$-space.  (In this case, Tits says that $y$ and $x \# V$ are ``li\'es'' in \cite[3.9]{Ti:Rsp}.)  For example, the vector $y$ from Example \ref{crucial.eg2} and ${\hw} \# V$ are connected because $y$ and $\hw$ are in the $\alpha_2$-space $V_{\alpha_2}$.  In contrast, the lowest weight vector $\lw$ and ${\hw} \# V$ are not connected, as $\phi(-\omega_6 + \mu)$ is a weight of $V$ for every weight $\mu$ of ${\hw} \# V$.  We find:

\begin{lem}[Cf.~\protect{\cite[3.16]{SpV}}] \label{dim.eq}
For $x$ and $y$ nonzero weight vectors in $V$,
\[
\dim \{ x, y, V \} \begin{cases} 
= 0 & \text{if $y$ and $x \# V$ are incident} \\
= 6 & \text{if $y$ and $x \# V$ are not incident but are connected} \\
\ge 17 & \text{if $y$ and $x \# V$ are neither incident nor connected.}
\end{cases}
\]
When the dimension is $6$, $\{ x, y, V \}$ is an $\alpha_2$-space.
\end{lem}

\begin{proof}
First consider the case where $x$ is the highest weight vector $\hw$.  We remind the reader that $y$ and ${\hw} \# V$ are incident if and only if $y$ is contained in ${\hw} \# V$, so the first equality follows from \ref{Va1}.  By Lemma \ref{orbits}, we are reduced to considering the cases where $y$ has weight $\la_2$ or $-\omega_6$.  
The second equality and the inequality now follow by Examples \ref{crucial.eg2} and \ref{crucial.eg1}, respectively.

In the general case, since the Weyl group acts transitively on the weights of $V$, there is some $n \in G$ such that $nx$ is a nonzero multiple of $\hw$. The dimension of $\{ x, y, V \}$ and the properties of $y$ and $x \# V$ being incident or connected are not affected by replacing $x$ and $y$ with $nx$ and $\phi(n) y$, respectively.  The general case now follows from the previous paragraph.
\end{proof}

\begin{rmk}
Since $G$ acts ``strongly transitively'' on the geometry $\G_V$ (see \cite[3.2.6]{Ti:BN}), it is sufficient to only consider the case where $x$ and $y$ are weight vectors.  That is, \eqref{dim.eq} holds for \emph{every} pair of singular vectors $x$, $y$.   We will not use this fact.
\end{rmk}

\begin{lem} \label{awkward}
Fix a proper inner ideal $X$ of $V$ such that $\dim X \ne 6$.  If $X$ has a basis $\B$ consisting of weight vectors, then 
   \[
   \psi(X) = \bigcap\nolimits_{b \in \B} b \# V
   \]
and $\{ X, \psi(X), V \} = 0$.
\end{lem}

\begin{proof}
Since $X$ has a basis consisting of weight vectors, so does $\psi(X)$.  Fix a $b \in \B$ and a weight vector $y \in \psi(X)$.  We claim that $\{ b, y, V \}$ is the zero subspace.  Otherwise, by Lemma \ref{dim.eq},  $\{ b, y, V \}$ is an $\alpha_2$-space or has dimension at least 17.  In particular, the dimension of $X$ is at least 7.  This implies that $X$ is a hyperline, so it cannot contain an $\alpha_2$-space by Remark \ref{hyper.quad}, and we have a contradiction.  This proves that $\{ b, y, V \}$ is zero.
Letting $y$ vary, we find that $\{ b, \psi(X), V \}$ is zero, which proves the second equation.  By Lemma \ref{dim1},
 $\psi(X)$ is contained in $b \# V$.

Conversely, suppose that $y$ is in the intersection of the $b \# V$'s.  Then $\{ x, y, V \}$ is zero for every $x \in X$.  Since the $b$'s span $X$, $y$ is in $\psi(X)$.  This proves the displayed equation.
\end{proof}

We can now compute $\psi(V_{\alpha_i})$ for $i \ne 1, 2$.

\smallskip

\begin{borel}{Computation of $V_{\alpha_3}$ and $V_{\alpha_4}$} \label{Va34}
The space $V_{\alpha_3}$ is spanned by the highest weight vector $\hw$ and a vector $x$ of weight $\la_3$.  We wish to compute $\psi(V_{\alpha_3})$, which is $V_{\alpha_6} \cap (x \# V)$ by Lemma \ref{awkward}. 
Each weight $\tau$ of $V_{\alpha_6}$ is a weight of $x \# V$ if and only if $\phi(\tau) - \la_3$ is a weight of $V$.  The five weights $\tau$ of $V_{\alpha_6}$ with a 1 as their last coordinate cannot belong to $x \# V$ because $\phi(\tau) - \la_3$ has a 2 as its first coordinate.  That is, $\psi(V_{\alpha_3})$ is contained in $V_{\alpha_5}$.

Since $\phi(\la_5) - \la_3 = -\phi(\la_2)$ is a weight of $V$, the weight $\la_5$ belongs to $x \# V$.  Figure \ref{weights} shows that $-\phi(\la_2) + \alpha_2$ is also a weight of $V$, hence $\la_5 + \alpha_2$ belongs to $x \# V$.  Continuing in this manner, we find that $V_{\alpha_5}$ is contained in $x \# V$, hence that $\psi(V_{\alpha_3})$ is $V_{\alpha_5}$.

The space $V_{\alpha_4}$ is spanned by $V_{\alpha_3}$ and a vector $y$ of weight $\la_4$.  The two weights $\tau$ of $V_{\alpha_5}$ that do not belong to $V_{\alpha_4}$ each have a 1 as their 5th coordinate, hence $\phi(\tau) - \la_4$ has a 2 as its 3rd coordinate, and such a $\tau$ does not belong to $\psi(V_{\alpha_4})$.  The three weights of $V_{\alpha_4}$ are easily checked to be weights of $y \# V$, hence $\psi(V_{\alpha_4})$ is $V_{\alpha_4}$.
\end{borel}

\begin{borel}{Computation of $V_{\alpha_5}$ and $V_{\alpha_6}$} \label{Va56}
By Lemma \ref{awkward} and \ref{Va34}, $\psi(V_{\alpha_5})$ is contained in $V_{\alpha_4}$.  Moreover, the 5th coordinate of $\phi(\la_4) - \la_5$ is $-2$, hence $\psi(V_{\alpha_5})$ is contained in $V_{\alpha_3}$.

Equation \eqref{lambdas} gives that the weight $\la_3$ belongs to $z \# V$ for $z$ a nonzero vector of weight $\la_5$.  Also,
\[
\la_1 = \la_3 + \alpha_1 = \phi(\la_5 + (\la_2 + \alpha_6)),
\]
so $\hw$ belongs to $z \# V$.  Similar calculations show that $V_{\alpha_3}$ is contained in $x \# V$ for $x$ of weight $\la_5 + \alpha_2$, hence $V_{\alpha_3}$ is equal to $\psi(V_{\alpha_5})$.

Lemma \ref{awkward} gives that $\psi(V_{\alpha_6})$ is contained in the 2-dimensional space $\psi(V_{\alpha_5}) = V_{\alpha_3}$.  The 6th coordinate of $\phi(\la_3) - \la_6$ is $-2$, hence $\la_3$ does not belong to $u \# V$ for every vector $u$ of weight $\la_6$, and $\psi(V_{\alpha_6})$ is contained in the $k$-span of the highest weight vector $\hw$.

We now show that $\hw$ is in $\psi(V_{\alpha_6})$.  Linearizing the identity $x^{\#\#} = N(x)x$ as in \cite{McC:FST} again (and going through McCrimmon's Equation (19)), we find the identity
\begin{equation}
z \# (y \# (x \# z)) = b(x, z^\#)y + b(x,y)z^\# + b(y,z)(x \# z) - x \# (y \# z^\#),
\end{equation}
which holds for every $x, y, z \in V$.  Substituting $z \mapsto \hw$, we find
   \[
   {\hw} \# (y \# (x \# {\hw})) = b(y,\hw) (x\# {\hw}).
   \]
Recalling that $b({\hw} \# x, \hw)$ is zero, we find that $\{ {\hw} \# x, \hw, y \}$ is zero for all $x, y \in V$.  That is, $\hw$ is in $\psi(V_{\alpha_6})$ and $\psi(V_{\alpha_6})$ is $V_{\alpha_1}$.
\end{borel}

We have proved that $\psi(V_{\alpha_i}) = V_{\phi(\alpha_i)}$ for all $i$.  In particular, the image of the standard chamber under $\psi$ is just the standard chamber.  This proves that $\psi$ is the automorphism of $\G_V$ induced by the automorphism $\phi$ of $G$.

In the language of classical projective geometry, $\psi$ is a \emph{hermitian polarity}.  Indeed,
since $\phi^2$ is the identity on $G$, $\psi^2$ is the identity on $\G_V$, i.e., $\psi$ is a polarity.  One says that $\psi$ is hermitian because the ``point'' $V_{\alpha_1}$ is contained in its ``polar'' $V_{\alpha_6}$.

\smallskip

It is natural to ask what $\psi$ does to inner ideals that are not  in $\G_V$, namely the 4-dimensional singular and 5-dimensional non-maximal singular subspaces.

\begin{eg} \label{closure}
In \cite[2.3]{Ti:Rsp}, Tits defines a geometry called an ``R-espace''.  It contains our geometry $\G_V$ plus two other types of objects: the 4- and 5-dimensional subspaces found as intersections in \ref{incidence.E}.  (Algebraically, the objects of Tits's R-espace are the nonzero, proper inner ideals in an Albert algebra.)  We now compute $\psi$ on these ``extra'' subspaces.

Let $X$ be the intersection of $V_{\alpha_2}$ and $V_{\alpha_5}$.  It is spanned by $V_{\alpha_4}$ and a nonzero vector $x$ of weight $(0, 1, 0, -1, 1, 0)$.
By Lemma \ref{awkward}, we have:
   \[
   \psi(X)  = \psi(V_{\alpha_4}) \cap (x \# V) = V_{\alpha_3}
   \]
Note that the 4-dimensional $X$ is properly contained in the 5-dimensional $\psi(\psi(X)) = V_{\alpha_5}$.  These results are what one expects from Tits's perspective \cite[3.3]{Ti:Rsp}.  But algebraically they are in contrast with $V_\delta$ for $\delta \in \D$, for which we have $\psi(\psi(V_\delta)) = V_\delta$.

The subspace $Y$ from \ref{incidence.E}---spanned by $X$ and a vector with weight $(0, 1, 0, 0, -1, 1)$---is similar.  The analogous computation shows that $\psi(Y)$ is $V_{\alpha_1}$, hence $\psi(\psi(Y))$ is the hyperline $V_{\alpha_6}$.
\end{eg}

Combining Lemma \ref{awkward} with the fact that every object in $\G_V$ is in the $G$-orbit of $V_\delta$ for some $\delta$, we obtain the following corollary:

\begin{cor} \label{not6}
If $X \in \G_V$ is not of dimension $6$, then 
   \[
   \psi(X) = \{ z \in V \mid \{ X, z, V \} = 0 \}. 
   \]
\end{cor}

Faulkner discussed the geometry $\G_V$ in terms of the brace product in \cite{Faulk:oct}, although he focussed on the points ($\alpha_1$-spaces) and hyperlines.  He described the duality on points and hyperlines by the equation displayed in the corollary.  However, that definition does not work for our purposes.  To see this, note that the 
 the set
   \[
   \{ z \in V \mid \{V_{\alpha_2}, z, V\} = 0 \}
   \]
is the zero space by computations as in Example \ref{closure}.
But the polar of $V_{\alpha_2}$ must be a 6-dimensional inner ideal. 


\begin{prop} \label{LN.prop}
For $X \in \G_V$, we have
   \[
   b(X, \psi(X)) = 0 \quad \text{and} \quad \{ \psi(X), X, \psi(X) \} = 0.
   \]
\end{prop}

\begin{proof}
By the transitivity of the $G$-action, we may assume that $X$ is $V_{\alpha_i}$ for some $i$.  Let $j$ be such that $\alpha_j = \phi(\alpha_i)$.  Further, let $d_i$ be such that $\omega_1 - d_i = \la_i$; it is a sum of positive roots.

\smallskip

We first argue that $b(X, \psi(X))$ is zero.  By \eqref{b.prop}, it can only be nonzero if there are vectors $x \in X$ and $x' \in \psi(X)$ of weights $\la$ and $\la'$ such that $\la + \phi(\la') = 0$.  But every weight of $X$  (resp.~$\psi(X)$) is at least $\la_i$ (resp.~$\la_j$), and
  \[
  \la_i + \phi(\la_j) = (\omega_1 + \omega_6) - (d_i + \phi(d_j));
  \]
we will show that this is $> 0$ for all $i$.
Consulting the tables in \cite{Bou:g4}, we find:
  \[
  \omega_1 + \omega_6 = 2\alpha_1 + 2\alpha_2 + 3\alpha_3 + 4\alpha_4 + 3 \alpha_5 + 2 \alpha_6.
  \]

When $i = 1$, $d_1$ is zero and $d_6$ is a sum of positive roots with no root occurring more than twice, as can be seen from Figure \ref{weights}.  Therefore, $\la_1 + \phi(\la_6) > 0$.  Applying $\phi$ to both sides of the equation, we find that $\la_6 + \phi(\la_1) > 0$.
When $i = 2, 3, 4$, or 5, no root appears more than once in $d_i$.  But $2 \le j \le 5$, hence the same is also true of $d_j$.  Therefore no root appears in $d_i + \phi(d_j)$ more than twice.  We have proved that $\la_i + \phi(\la_j) > 0$ for all $i$, hence $b(X, \psi(X))$ is necessarily zero.

\smallskip

We now prove that the second equation holds.  The space $\{ \psi(X), X, \psi(X) \}$ is a direct sum of its weight spaces, and each weight $\mu$ is at least $\phi(\la_i) + 2 \la_j$ by \eqref{duality.br}.  That is, each weight $\mu$ is of the form $\omega_1 - d$ for 
\begin{equation} \label{disp.1}
0 \le d \le \omega_1 - (\phi(\la_i) + 2 \la_j) = (\phi(d_i) + 2d_j)-(\omega_1 + \omega_6).
\end{equation}

It follows from the definition of $\la_j$ as the lowest weight of $V_{\alpha_j}$ that $d_j$ does not include the root $\alpha_j$.  Similarly, $\phi(d_i)$ also does not include $\alpha_j$.  Therefore, the coefficient of $\alpha_j$ on the right side of \eqref{disp.1} is negative.  In particular, the equation $d \ge 0$ is impossible, and $\mu$ cannot be a weight of $V$.  This proves that $\{ \psi(X), X, \psi(X) \}$ is zero.
\end{proof}

 In \cite[p.~260]{LN}, Loos and Neher defined
  \[
  \inid(X) := \{ y \in V \mid \text{$\{ y, X, y \} = 0$ and $\{ X, y, V \} \subseteq X$} \},
  \]
for $X$ a subspace of $V$.  Clearly, $\psi(X)$ contains $\inid(X)$, and the preceding proposition shows that the two concepts agree for $X \in \G_V$.

\begin{rmk} \label{brace.rmk}
We used the brace product in this section so that we could apply results from Jordan theory; for an explanation of why Jordan theorists like it, see \cite[pp.~7ff]{McC}.  There is another ternary product that is perhaps more natural from the perspective of representation theory.  Up to a scalar multiple, there is a unique $G$-invariant linear map $\langle \rangle$ from $V \ot V^*$ to the Lie algebra of $G$.   The composition
\[
[ x, y, z ] := {\langle x, h(y) \rangle}\ z
\]
satisfies an equation analogous to \eqref{duality.br}.  Moreover, if one choose the appropriate multiple of $\langle \rangle$, then
\[
[x , y ,z ] = \{ x, y, z \} - \frac23 b(x,y) z \quad \text{for $x, y, z \in V$.}
\]
(This product was studied in \cite[26.9ff]{Frd:E7.8} and \cite{SpV}.)
A subspace of $V$ is an inner ideal with respect to the brace product if and only if it is an inner ideal with respect to the bracket product.  If we define $[\psi](X)$ to be the set of all $y \in V$ such that $[X, y, V] \subseteq X$, then $\psi(X)$ is contained in $[\psi](X)$ for $X \in \G_V$ by Prop.~\ref{LN.prop}.  In fact, $\psi$ and $[\psi]$ agree for $\alpha_i$-spaces with $i \ne 2$ by \cite[5.4, 6.7]{G:struct}.

We remark that the map $z \mapsto \{ x, y, z \}$ belongs to the Lie algebra of the \emph{structure group} of $V$.  That group is a reductive envelope of $G$ with center a rank 1 torus.
\end{rmk}

\begin{rmk} \label{jord.free}
A Jordan theorist would start not with a split group of type $E_6$ but rather with a ``cubic norm structure'' (a characteristic-free version of an Albert algebra) consisting of the 27-dimensional vector space $V$, the map $x \mapsto x^\#$, the cubic form $N \!: V \ra k$, and a base point $1 \in V$, all satisfying certain axioms \cite[\S1]{McC:FST}.  From that perspective the results of this paper hold without the restriction that $k$ has characteristic zero.  First, the subgroup of $GL(V)$ of elements that preserve $N$ form a  split simply connected group of type $E_6$ \cite[p.~151]{Sp:jord}.  Sections \ref{Tits} through \ref{realize.2} go through with no change, see Remark \ref{min.free}.  For \S\ref{E6}, there is no need to construct $N$ and $\#$.  The form $b$ is obtained as a logarithmic derivative of $N$, see e.g.\ \cite[1.18]{Sp:jord}.  The material from \ref{singular} through the end of \S\ref{E6} goes through with only cosmetic changes.  All of \S\ref{E6.dual} goes through when the characteristic of $k$ is different from 2.
(To get theorems that work in all characteristics, one should consider not the brace product but rather a quadratic map whose bilinearization is the brace product.  We have not done so here in order to make the exposition somewhat smoother.)
\end{rmk}

\medskip

{\small
\noindent\textbf{Acknowledgements.} We are indebted to H.P.~Petersson for pointing out the references \cite{LN} and \cite{Racine}, to D.Z.~Djokovic and J-P.~Serre for their comments, and to E.~Neher for his advice about Jordan algebras.}




\providecommand{\bysame}{\leavevmode\hbox to3em{\hrulefill}\thinspace}
\providecommand{\MR}{\relax\ifhmode\unskip\space\fi MR }
\providecommand{\MRhref}[2]{%
  \href{http://www.ams.org/mathscinet-getitem?mr=#1}{#2}
}
\providecommand{\href}[2]{#2}

\end{document}